\documentclass[letterpaper,10pt,twocolumn]{article}
\usepackage{cite}
\usepackage[margin=1in]{geometry}
\usepackage{amsmath,amssymb,amsfonts}
\usepackage{enumitem}
\usepackage{algorithm,algorithmic}
\usepackage{textcomp}
\usepackage{verbatim}
\usepackage{xcolor}
\usepackage{mathtools}

\newtheorem{assumption}{Assumption}
\newtheorem{theorem}{Theorem}
\newtheorem{lemma}[theorem]{Lemma}
\newtheorem{prop}[theorem]{Proposition}
\newtheorem{corollary}[theorem]{Corollary}

\newtheorem{definition}{Definition}
\newtheorem{problem}{Problem}
\newtheorem{remark}{Remark}

\newcommand{\RR}{\mathbb{R}}
\newcommand{\CC}{\mathbb{C}}

\newcommand{\cB}{\mathcal{B}}
\newcommand{\cC}{\mathcal{C}}
\newcommand{\cCp}{\mathcal{C}^{\perp}}

\newcommand{\cK}{\mathcal{K}}
\newcommand{\cL}{\mathcal{L}}

\newcommand{\cP}{\mathcal{P}}

\newcommand{\cT}{\mathcal{T}}

\newcommand{\cV}{\mathcal{V}}
\newcommand{\cW}{\mathcal{W}}

\newcommand{\al}{\alpha}

\newcommand{\ale}{{\alpha}_e}
\newcommand{\be}{{\beta}_e}

\newcommand{\eps}{\epsilon}
\newcommand{\stol}{\sigma_{\text{tol}}}
\newcommand{\shattol}{\hat{\sigma}_{\text{tol}}}
\newcommand{\veps}{\varepsilon}

\newcommand{\cs}{c_{\star}}
\newcommand{\Fb}{\bar{F}}
\newcommand{\Fbz}{\bar{F_{\circ}}}
\newcommand{\Gbz}{\bar{G_{\circ}}}
\newcommand{\Gb}{\bar{G}}

\newcommand{\nh}{\hat{n}}
\newcommand{\nb}{\bar{n}}

\newcommand{\Sigh}{\hat{\Sigma}}
\newcommand{\Sigt}{\tilde{\Sigma}}
\newcommand{\thh}{\hat{\theta}}

\newcommand{\tht}{\tilde{\theta}}

\newcommand{\xih}{\hat{\xi}}

\newcommand{\xipe}{{\xi}_{pe}}
\newcommand{\xis}{{\xi}_{\star}}
\newcommand{\xisone}{{\xi}_{\star1}}
\newcommand{\xistwo}{{\xi}_{\star2}}
\newcommand{\xit}{\tilde{\xi}}
\newcommand{\yb}{\bar{y}}

\newcommand{\ol}{\overline}
\newcommand{\pe}{_{pe}}

\newcommand{\pr}{_{\perp}}
\newcommand{\inv}{^{-1}}
\newcommand{\0}{_{\circ}}
\newcommand{\x}{\times}
\newcommand{\T}{^{\intercal}}
\newcommand{\p}{\partial}

\newcommand{\AVG}[1]{\lim_{T \to \infty} \frac{1}{T} \int_{t_0}^{t_0 + T} #1 \,d\tau}

\newcommand{\tqed}{\hfill$\triangleleft$}

\newcommand{\red}{\color{red}}

\newcommand{\blue}{\color{blue}}

\DeclareMathOperator{\Span}{span}

\DeclareMathOperator{\Img}{Im}

\DeclareMathOperator{\Ker}{Ker}

\DeclareMathOperator{\col}{col}
\DeclareMathOperator{\diag}{diag}

\DeclareMathOperator{\rnk}{rank}

\DeclareMathOperator{\tr}{trace}
\DeclareMathOperator{\vect}{vec}

\DeclareMathOperator{\eig}{eig}
\DeclareMathOperator{\roots}{roots}
\usepackage{subcaption}
\usepackage{graphicx}

\begin{document}

\title{Continuous Time Identification of Linear Systems: Extended Version}
\author{
Fatima J. Ghadieh$^{1}$ \and
Mireille E. Broucke$^{1,2}$
\thanks{$^{1}$Supported by the Natural Sciences and Engineering Research Council of Canada (NSERC). Electrical and Computer Engineering, University of Toronto, Toronto, ON, Canada.}%
\thanks{$^{2}$Email: \texttt{broucke@control.utoronto.ca}}%
}
\maketitle

\begin{abstract}
We consider a problem of developing a framework for model identification adhering to the tenets of neuromorphic computation, without resorting to neural networks as the mathematical substrate. In particular, all computations take place in continuous time. We are naturally led to adaptive observers, where the main technical obstacle is the possible mismatch between the unknown plant order and the observer order. The key concept that informs the proposed framework is an overparameterized model, an input-output equivalent model that provides a suitable parameterization in the overmodeled case, with theoretical extensions also addressing the undermodeled case. A discrete algorithm orchestrates successive experiments to incrementally learn the model order, while a standard parameter adaptation law learns the parameters.
\end{abstract}

\section{Introduction}
\label{sec:introduction}

This paper initiates an investigation into a new area of control theory called {\em neuromorphic system identification}, which regards how the brain learns models. While there is a plethora of research on learning using neural networks, there has been no dedicated effort in control theory to formalize model learning in the brain at a systems level, without neural networks as the mathematical substrate \cite{VIDYASAGAR97}. Such neuroscience applications have not been well-served by current system identification frameworks, which assume computations suitable for a digital computer. This clear research gap opens a door for a fresh perspective. 

We consider accepted tenets of neuromorphic computation that capture the capabilities of neural networks without explicitly modeling them:
(i) 
Signal processing in the brain occurs in continuous time by frequency modulated nerve impulses. Only specialized motor systems such as the saccadic and visuomotor systems are capable of data sampling. 
(ii) 
Recursive computations dominate over storage and manipulation of large data sets. 
(iii)
Both individual neurons as well as neural networks are capable of performing covariance and singular value decomposition (SVD) computations recursively \cite{OJA82,BALDI89,YANHELMKEMOORE94}.
(iv)
The brain manages noise through massively parallel signal processing \cite{TABAREAU10}. 
(v)
Neural networks are capable of implementing the gradient law (or Widrow-Hoff rule) \cite{DAYAN01} and can operate as adaptive filters \cite{FUJITA82,DEAN10}.
The central research challenge is to identify a control-theoretic system identification framework that meets some, if not all, of the tenets of neuromorphic computation. 

\subsection{Literature Review}

The modern era of system identification was initiated with the Ho-Kalman algorithm for discrete time systems \cite{HOKALMAN66}. Major developments followed in the 1980s-1990s including subspace identification methods, instrumental variable methods, and frequency domain methods  \cite{SODERSTROM89,VANOVERSCHEE96book,LJUNG99,VERHAEGEN07,KATAYAMA05,PINTELON12}.
A resurgence of interest has been fueled by {\em data-driven methods} \cite{VANWAARDE25}, where the concern is with the quality and amount of data collected \cite{CARE18,SARKAR21}. 
These methods primarily regard discrete time systems, whereas here we study continuous processes. Direct identification of continuous time systems has been widely studied \cite{JOHANSSON99,GARNIER08,AFRI15,GARNIER15,GONZALES21}, with notable recent progress in \cite{HUANG22,HUANGFENG22}, resulting in methods that use discretely sampled input-output data to construct continuous time models. Unfortunately, we are not able to directly adopt these methods since our neuromorphic framework disallows data sampling.

System identification of continuous time systems using purely continuous time signal processing naturally leads to methods involving {\em adaptive observers} \cite{KREISSEL77,NARENDRA89}, which have been of sustained research interest for many decades \cite{AFRI17,FARZA09,PYRKIN23,TOMEI23,TYUKIN13}, 
including in neuroscience applications \cite{BURGHI24}. The key obstacle in applying adaptive observers for system identification is that the plant order must be known, otherwise numerical difficulties arise \cite{LAI24}.

\begin{figure*}[!t]
\centering
\includegraphics[width=2\columnwidth]{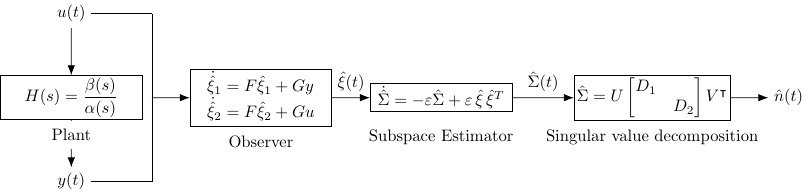}
\caption{Stages of identification: continuous time measurement of input-output data from the plant, two filtering stages, and an SVD computation to yield an estimate of the model order. A parameter identification stage is not depicted.}
\label{fig:blockdiagram}
\end{figure*}

The main inspiration for the methods of this paper is \cite{KREISSEL96}. They put forward the idea of an overparameterized model, a gateway to converting system identification to a problem of parameter adaptation, by providing a suitable parameterization in the overmodeled case. The only comparable method that we are aware of for working directly with models of mismatched order is \cite{LAI24}, where system identification of discrete time systems is addressed using recursive least squares (RLS), with the main challenge being a lack of persistent excitation of a regressor. 

\subsection{Contributions}

After performing a broad literature survey, we found we could not directly adopt existing methods of system identification, whereas continuous time adaptive observers do provide a plausible framework for neuroscience applications. In this setting, the main challenge is to identify the plant order, since recovery of parameters for a SISO LTI plant of known order is a solved problem \cite{KREISSEL77,NARENDRA89}.
The key contribution of the paper is a framework, seen in Figure~\ref{fig:blockdiagram}, which consists of three stages of processing to recover the plant order and an additional stage (not shown) to estimate plant parameters. While the framework utilizes familiar constructs such as observers and SVD computations, the key distinction is that the observer and plant need not be matched in their orders. 

Section~\ref{sec:overp} draws attention to elegant but overlooked work in \cite{KREISSEL94,KREISSEL96} on overparameterized models. These are state space models whose order may exceed the plant order yet they generate the same input-output behavior and the same output response under suitable initial conditions, thus capturing the notion of a non-minimal realization in parametric form. Several remarkable properties of overparameterized models are presented in Section~\ref{sec:overp}, some of which were reported in \cite{KREISSEL96}, but we also present some new results. Section~\ref{sec:order} presents one of the main contributions of the paper,  Theorem~\ref{thm:order}, which characterizes the model order in terms of the dimension of the controllable subspace of the overparameterized model. 

In Section~\ref{sec:contract} we use a Kreisselmeier observer to recover a special solution of the overparameterized model that lives in its controllable subspace. It is shown that with a sufficiently rich input, the observer state recovers the full controllable subspace of the overparameterized model. We apply a subspace estimator from \cite{UZEDA23,UZEDA25}, which was originally developed to estimate PE subspaces for robust adaptive control. A robust computation of the dimension of the controllable subspace is obtained using the SVD.

Section~\ref{sec:underparam} removes the requirement that an upper bound on the plant order is known. Our main result extends results on persistent excitation of the observer state to a novel case when the observer under-models the plant.
Section~\ref{sec:algorithm} presents a model order identification algorithm that incrementally estimates the plant order, a scheme that resembles the traditional discrete time Ho-Kalman algorithm. Section~\ref{sec:params} addresses the well-studied problem of parameter identification when the plant order is known, however from a different perspective, as we exploit the existence of the special solution of the overparameterized model to streamline arguments about PE of regressors. Section~\ref{sec:examples} discusses the robustness of the method and provides examples. A {\em dimension plot} aids in visualizing the model order estimation process, and several examples highlight the efficacy of the method. 
 
\section{Preliminaries}
\label{sec:prelim}

We use the following standard notation. 
The disjoint union of two sets is denoted $\sqcup$. 
For a subspace $\cV \subset \RR^n$, its dimension is denoted $\dim(\cV)$. 
For a matrix $M$, its image and kernel are denoted $\Img(M)$ and  $\Ker(M)$. For a matrix $A \in \RR^{n \times n}$, its spectrum, rank, and transpose are denoted as $\eig(A)$, $\rnk(A)$, and $A\T$. The $i$th largest eigenvalue and singular value of a symmetric matrix $A \in \RR^{n \times n}$ are denoted as $\lambda_i(A)$ and $\sigma_i(A)$. The degree of polynomial $f$ is denoted as $\deg(f)$. If $u(t)$ is a function of time, then $u(s)$ denotes its Laplace transform. 
Let $\cB(r) \subset \RR^n$ denote the open ball of radius $r \ge 0$, and  let $\| x \|_{\cL_\infty}$ denote the $\cL_\infty$ norm of a signal $x(t) \in \RR^n$ (assuming it exists).
We use the following canonical matrices (whose dimensions depend on the context):

\begin{equation}
\begin{aligned}\label{eq:ABC2}
    A\0 &= 
    \left[\begin{array}{cccc}
        0 & \multicolumn{3}{|c}{} \\
        \vdots & \multicolumn{3}{|c}{I} \\
        0 & \multicolumn{3}{|c}{} \\ \hline
        0 & \multicolumn{2}{c}{\cdots} & 0
    \end{array}\right] \,,
    \
    B\0 =
    \left[\begin{array}{c}
        0 \\ \vdots \\ 0 \\ 1
    \end{array}\right] \\
    C\0 &=
    \left[\begin{array}{cccc}
        1 & 0 & \cdots & 0
    \end{array}\right] \,.
\end{aligned}
\end{equation}

We review some background from \cite{NARENDRA89,SASTRY89}. 
A bounded signal $w(t) \in \RR^q$ is {\em persistently exciting} (PE) if there exist $\beta\pe$ and $T\pe > 0$ such that
\begin{align}\label{eq:PE}
\frac{1}{T\pe} \int_{t}^{t+T\pe} w(\tau) w\T(\tau) \,d\tau  
\succeq \beta\pe I \,, \quad \forall \ t \geq 0 \,.
\end{align}
The pair $(\beta\pe,T\pe)$ is called the {\em PE constants of $w$} (they are not unique). We will generally be concerned with the following class of signals.

\begin{definition} \label{def:auto}
A signal $w(t) \in \RR^q$ is said to be {\em stationary} if it is bounded and piecewise continuous, and its {\em autocovariance matrix}

\begin{align*}
R_w( \rho ) := \AVG{w(\tau) w\T(\tau + \rho )}
\end{align*}
exists, is independent of the initial time $t_0$, and the convergence of the time average is uniform in $t_0 \geq 0$. 
\tqed
\end{definition} 
In the sequel we only require the zero lag autocovariance, denoted as $R[w] = R_w(0)$.  If a signal $w:\RR \rightarrow \RR$ has a sinusoidal component at frequency $\omega_0$ then its {\em spectral measure} $S_w(\omega)$ has \textit{point masses} at $\omega_0$ and $-\omega_0$. A stationary signal $w(t)$ is \textit{sufficiently rich} (SR) of order $k$ if $S_w(\omega)$ contains at least $k$ points \cite[p.~255]{IOANNOU12}.

\subsection{PE Subspaces}
\label{eq:PEsubspace}

Consider a stationary signal $w(t) \in \RR^q$.
The {\em PE subspace} of $w$ is defined to be $\cW = \Img ( R[w] )$, while the {\em non-PE subspace} of $w$ is $\cW^{\perp}$. 
We denote $q\pe := \dim(\cW)$ as its {\em degree of persistent excitation}. 
When $q\pe = q$, $w$ is PE by \cite[Proposition~2.7.1]{SASTRY89}, whereas when $q\pe = 0$, $w$ has {\em no persistent excitation}. 
Stationary signals are amenable to a PE decomposition.
\begin{prop}[Prop~1,\cite{UZEDA25}]
\label{prop:PEdecomp}
Suppose $w$ is stationary. If $1 \leq q\pe < q$, let $\begin{bmatrix}W & W\pr\end{bmatrix} \in \RR^{q \x q}$ be orthogonal such that 
$\cW = \Img(W)$ and $\cW^{\perp} = \Img(W\pr)$.
Then the {\em PE decomposition of $w$} is
\begin{align}
\label{eq:PEdecomp}
w &= W W\T w + W\pr W\pr\T w =: W w\pe + W\pr w\pr \,,
\end{align}
where $w\pe(t) \in \RR^{q\pe}$ is PE and $w\pr(t) \in \RR^{(q - q\pe)}$ has no persistent excitation.
\tqed
\end{prop}

The conceptual underpinning of our system identification method is that the controllable subspace of an overparameterized model of the plant contains all information about the plant order and parameters. To gain access to that controllable subspace, we use a sufficiently rich input to ensure that the PE subspace of the state of the overparameterized model is equivalent to its controllable subspace. We can then recover the controllable subspace by using subspace estimation tools for PE subspaces from \cite{UZEDA25}.  
Consider the LTI system
\begin{subequations}
\label{eq:ABC}
\begin{align}
\label{eq:AB}
\dot{x} &= A x + B u \\
y       &= C x \,,
\end{align}
\end{subequations}
where $x(t) \in \RR^n$, $u(t) \in \RR$, and $y(t) \in \RR$. The key notion is that for LTI systems there is an intrinsic link between PE subspaces and controllable subspaces.
\begin{lemma}
\label{lem:PEC}
Consider \eqref{eq:AB} with $A$ Hurwitz. Suppose $u$ is stationary and SR of order $n$. Then the PE subspace of $x$ is 
$\cC(A,B) = \Img( \begin{bmatrix} B & AB & \cdots & A^{n-1} B \end{bmatrix})$, the controllable subspace of \eqref{eq:AB}.
\tqed
\end{lemma}

\section{Problem Statement}

Consider a SISO LTI system modeled as a transfer function
\begin{align}
\label{eq:theplant}
H(s) &= \frac{\beta(s)}{\al(s)} =: 
\frac{\beta_{n-1} s^{n-1} + \ldots + \beta_1 s + \beta_0} 
     {s^n + \al_{n-1} s^{n-1}  + \ldots + \al_1 s + \al_0} \,,
\end{align}
where $n$ is the order of the plant. 
The following standing assumptions hold throughout the paper.
\begin{assumption}
\label{assum:plant}
The plant \eqref{eq:theplant} satisfies the following.
\begin{enumerate}[label=(A\arabic*), leftmargin=*, ref=(A\arabic*)]
\item \label{as:coprime}
$\{ \al(s), \beta(s) \}$ are coprime.
\item
 $\al(s)$ is Hurwitz.
\item
The model order $n$ and parameters $\al := ( \al_0, \ldots, \al_{n-1} ) \in \RR^n$ and $\beta := ( \beta_0, \ldots, \beta_{n-1} ) \in \RR^n$ are unknown.
\item
The measurements are input $u(t) \in \RR$ and output $y(t) \in \RR$.
\tqed
\end{enumerate}
\end{assumption}
\begin{problem}
Consider the system \eqref{eq:theplant} satisfying Assumption~\ref{assum:plant}. Develop a method  to estimate the plant order $n$ and the plant parameters using continuous input-output data.
\tqed
\end{problem}

\section{Overmodeled Plant}\label{sec:overp}

We begin our study with the overmodeled case.
\begin{assumption} \label{assum:over} 
An upper bound on the plant order $\nb \ge n$ is known.
Let $p := \nb - n$ be the excess order. 
\tqed
\end{assumption}

Consider the {\em overparameterized model} studied in \cite{KREISSEL96}:
\begin{subequations}
\label{eq:overp1}
\begin{align}
\dot{\xi}_1 &= F \xi_1 + G y\0 \\
\dot{\xi}_2 &= F \xi_2 + G u \\
y\0 & = \theta_1\T \xi_1 + \theta_2\T \xi_2 \,, 
\end{align}
\end{subequations}
where $(\xi_1, \xi_2) \in \RR^{2 \nb}$ is the state, $y\0 \in \RR$ is the overparameterized model output, and $\theta_1, \theta_2 \in \RR^{\nb}$ are to-be-determined parameters. We assume $F$ is Hurwitz, and $(F,G)$ is a controllable pair. Specifically we select
\begin{align}
\label{eq:FG}
F &= A\0 + B\0 (-f)\T \,, \qquad G = B\0 \,, 
\end{align}
where $f = (f_0, f_1, \ldots, f_{\nb-1}) \in \RR^{\nb}$ and 
$f(s) = s^{\nb} + f_{\nb-1} s^{\nb-1} + \cdots + f_1 s + f_0$ is a Hurwitz polynomial. Defining $\xi := (\xi_1,\xi_2) \in \RR^{2 \nb}$ and parameters $\theta = (\theta_1,\theta_2) \in \RR^{2 \nb}$, \eqref{eq:overp1} can be written more compactly as 
\begin{subequations}
\label{eq:overp2}
\begin{align}
\dot{\xi} 
&= \begin{bmatrix} (F + G \theta_1\T) & G \theta_2\T \\ 0 & F \end{bmatrix} \xi +
\begin{bmatrix} 0 \\ G \end{bmatrix} u =: \Fb(\theta) \xi + \Gb u \\
y\0 
&= \theta\T \xi \,.
\end{align}
\end{subequations}

We define the set of parameters for which the overparameterized model is a model of the plant. That is,
\begin{align}
\cP &:= 
\left\{ \theta \in \RR^{2 \nb} ~|~ 
\theta\T ( sI - \Fb(\theta) )\inv \Gb = \frac{\beta(s)}{\alpha(s)} \right\} \,.
\end{align}
The following result characterizes the elements of $\cP$.
\begin{prop}
\label{prop:thetas}
Suppose Assumptions~\ref{assum:plant} and \ref{assum:over} hold.
Consider \eqref{eq:overp1} with $(F,G)$ defined in \eqref{eq:FG}. Let $\delta(s)$ be any monic polynomial of degree $p = \nb - n$. Select $\theta = (\theta_1,\theta_2) \in \RR^{2 \nb}$ with
\begin{align}
\label{eq:thetas}
\theta_1 = f - \ale \,, \qquad \theta_2 = \be \,,
\end{align}
where $\ale = (\bar{\al}_0, \bar{\al}_1, \ldots, \bar{\al}_{\nb-1})$ and $\be = (\bar{\beta}_0, \bar{\beta}_1, \ldots, \bar{\beta}_{\nb-1})$ are such that
\begin{subequations}
\label{eq:extended}
\begin{align}
\nonumber
\ale(s) 
&:= s^{\nb} + \bar{\al}_{\nb-1} s^{\nb-1} + \cdots + \bar{\al}_1 s + \bar{\al}_0 \\
&= \al(s) \delta(s) \\ 
\be(s) 
&:= \bar{\beta}_{\nb-1} s^{\nb-1} + \cdots + \bar{\beta}_1 s + \bar{\beta}_0 = \beta(s) \delta(s) \,.
\end{align}
\end{subequations}
Then $\theta \in \cP$ if and only if $\theta$ satisfies \eqref{eq:thetas} for some $\delta(s)$.
\tqed
\end{prop}
Proposition~\ref{prop:thetas} captures the intuitive idea that for a higher-order SISO LTI model to match the input-output behavior of a lower order plant model, there must be pole-zero cancellations. We can go further to match the output transient response of the two models by choosing appropriate initial conditions of the overparameterized model. 

\begin{prop}
\label{prop:matching}
Suppose Assumptions~\ref{assum:plant} and \ref{assum:over} hold.
Consider $(F,G)$ defined in \eqref{eq:FG} and $\theta \in \cP$. Then for each initial condition $\yb(t_0) := (y,\dot{y},\ldots,y^{(n-1)})(t_0)$ of \eqref{eq:theplant}, there exists $\xi(t_0) \in \RR^{2 \nb}$ such that \eqref{eq:overp2} satisfies: 
\begin{itemize}
\item[(i)]
$y\0(t) = y(t)$ for all $t \ge t_0 \ge 0$. 
\item[(ii)]
$\xi(t) \in \cC(\theta)$, the controllable subspace of \eqref{eq:overp2}, for all $t \ge t_0 \ge 0$. \tqed
\end{itemize}
\end{prop}

The methodology to be developed involves using the controllable subspace $\cC(\theta)$ to extract information about the plant order and parameters. Unfortunately, there appears to be a chicken and egg problem, because the controllable subspace $\cC(\theta)$ depends on the choice of $\theta \in \cP$, but the elements of $\cP$ are unknown. A step toward resolving this dilemma is the observation that all members of $\cP$ give rise to the same controllable subspace. To that end, 
let $\xis(t;\theta,\yb(t_0))$ denote the output matching solution identified in Proposition~\ref{prop:matching}. 
The next result ensures that we may study the properties of the controllable subspace as well as the special solution $\xis(t;\theta,  \yb(t_0))$ of \eqref{eq:overp2} for any choice $\theta \in \cP$. 

\begin{prop}
\label{prop:controllable}
Suppose Assumptions~\ref{assum:plant} and \ref{assum:over} hold.
Consider \eqref{eq:overp2} with $(F,G)$ defined in \eqref{eq:FG}. For each initial condition $\yb(t_0) \in \RR^n$ of \eqref{eq:theplant} and for all $\theta, \theta' \in \cP$, we have
\begin{itemize}
\item[(i)]  
$\cC(\theta) = \cC(\theta')$.
\item[(ii)]   
$\xis(t;\theta,\yb(t_0)) = \xis(t;\theta',\yb(t_0))$ 
for all $t \ge t_0 \ge 0$. 
\end{itemize}
\tqed
\end{prop} 

We henceforth drop the dependence on $\theta$ and simply write $\cC$ to denote the controllable subspace for any choice of overparameterized model,
and $\xis(t;\yb(t_0))$ to denote the output matching solution.

\begin{remark}
The remarkable properties of the overparameterized model \eqref{eq:overp2} were discovered by Kreisselmeier and Lozano \cite{KREISSEL94,KREISSEL96}. Proposition~\ref{prop:thetas} is from \cite{KREISSEL96}; its proof has been included to make the paper self-contained. The proof of Proposition~\ref{prop:matching}, given in Section~\ref{proof:prop:matching}, uses a more streamlined state space argument compared to the Laplace domain arguments in \cite{KREISSEL96}. Proposition~\ref{prop:controllable} can be deduced by the results of \cite{KREISSEL96} but was not proven there.
\tqed
\end{remark}

\subsection{Plant Order Estimation}\label{sec:order}

The following is one of the main results of the paper.

\begin{theorem}
\label{thm:order}
Suppose Assumptions~\ref{assum:plant} and ~\ref{assum:over} hold.
Consider \eqref{eq:overp2} with $(F,G)$ defined in \eqref{eq:FG} and $\theta \in \cP$. The dimension of the controllable subspace of \eqref{eq:overp2} is 
$n_c := \dim(\cC) = n + \nb$.
\tqed
\end{theorem}

Theorem~\ref{thm:order} provides the grounds to determine the model order $n$. If we can estimate the dimension of the controllable subspace of \eqref{eq:overp2}, then we can recover $n$.
By combining Proposition~\ref{prop:controllable} with Theorem~\ref{thm:order} we obtain further insight on the eigenvalues of the overparameterized model.

\begin{corollary}
\label{cor:poles}
Suppose Assumptions~\ref{assum:plant} and ~\ref{assum:over} hold.
Consider \eqref{eq:overp2} with $(F,G)$ defined in \eqref{eq:FG} and $\theta \in \cP$. Then the controllable modes of \eqref{eq:overp2} are $\roots(\alpha(s)) \sqcup \roots(f(s))$, while the uncontrollable modes are $\roots(\delta(s))$.
\tqed
\end{corollary} 

\begin{remark}
Theorem~\ref{thm:order} belongs to a family of results on model order estimation, with the unique feature that it characterizes the model order in terms of a geometric property of an input-output equivalent model. As such, it places no requirements on the richness of the input or the quality or amount of data. By way of comparison, characterizations of the model order in terms of collected data are seen in Theorem~2 in \cite{MOONEN89} (see also Proposition~1.4 in \cite{VANWAARDE25});  Section~4 in \cite{HAVERKAMP96}; Theorem~1 in \cite{VANOVERSCHEE96}; Theorem~2 in \cite{VANOVERSCHEE96book}; Section~3.3 in \cite{BASTOGNE97}; and Theorem~5.4 of \cite{SARKAR21}.
\tqed
\end{remark}

\subsection{Overparameterized Observer}
\label{sec:contract}

There are two difficulties with applying Theorem~\ref{thm:order}. First, we do not have access to $\xis(t;\yb(t_0))$. Second, we do not have in hand a method to estimate the dimension of $\cC$, even if $\xis(t;\yb(t_0))$ were available. In this section we remove the first difficulty.
We consider an {\em overparameterized observer} of the form
\begin{subequations}
\label{eq:observer}
\begin{align}
\dot{\xih}_1 &= F \xih_1 + G y \\
\dot{\xih}_2 &= F \xih_2 + G u \,,
\end{align}
\end{subequations}
with $\xih := (\xih_1, \xih_2) \in \RR^{2 \nb}$ and arbitrary initial conditions $\xih(t_0) := (\xih_1(t_0), \xih_2(t_0))$. 
By Proposition~\ref{prop:matching}, $y\0(t) = y(t)$ for all $t \ge t_0 \ge 0$ with a corresponding solution $\xis(t;\yb(t_0))$. We use the following notation
\begin{align}
\label{eq:FboGbo}
\Fbz &:= 
\begin{bmatrix}
             F & 0\\
             0 & F\\
\end{bmatrix} ; \qquad 
\Gbz  := 
\begin{bmatrix}
             G & 0\\
             0 & G
\end{bmatrix}\,.
\end{align}
Define $\xit := \xih - \xis(t;\yb(t_0))$. Then we have
\begin{align}
\label{eq:xit}
\dot{\xit} &= \Fbz \xit  \,.
\end{align}
That is, the states of the overparameterized observer \eqref{eq:observer} converge exponentially to $\xis(t;\yb(t_0))$ of \eqref{eq:overp2}, even though the latter depends on unknown parameters and initial conditions. 

\begin{remark}
\label{rem:FGf}
The observer and corresponding overparameterized model were constructed using $(F,G)$ in controllable canonical form, a choice  informed by the ease of deriving \eqref{eq:thetas}, as well as simplifying algebraic manipulations in the proofs of Propositions~\ref{prop:thetas}-\ref{prop:controllable}. Other choices of $(F,G)$ are clearly possible, resulting in new formulas for  $\theta$. This flexibility raises an intriguing problem of whether one can identify canonical structures for $(F,G)$ that characterize, in an average sense, the massive parallel filtering of neural networks \cite{TABAREAU10}.
A further issue regards the selection of $f(s)$, which has received some attention in the system identification literature 
\cite{LJUNG02,GARNIER03,GARNIER08} in terms of optimal 
selection of pre-filters. Section~14.4 of \cite{LJUNG02} suggests that pre-filters can be interpreted as noise models, while \cite[Chapter~5]{GARNIER08} discusses the optimal selection of $f(s)$ in instrument variable methods. Whether these suggestions are relevant in brain modeling requires further investigation, so we have not further pursued the optimal selection of $f(s)$.
\tqed
\end{remark} 

\subsection{Subspace Estimation} 
\label{sec:subspace}

In this section we remove the second difficulty in applying Theorem~\ref{thm:order}; namely, we show how to gain access to the controllable subspace $\cC$ of the overparameterized model. 
To that end, let $T := \begin{bmatrix} U & U\pr \end{bmatrix} \in \RR^{2 \nb \times 2 \nb}$ be an orthogonal matrix such that
\begin{align} \label{eq:UUpr}
\cC = \Img(U), \qquad \cCp = \Img(U\pr) \,.
\end{align}
The main idea to recover $\cC$ is to inject sufficient excitation into the overparameterized model state via a sufficiently rich input so that the full controllable subspace is revealed through the time evolution of $\xis(t;\yb(t_0))$, based on the equivalence established in Lemma~\ref{lem:PEC}. 

\begin{prop}
\label{prop:PE}
Suppose $u$ is stationary and SR of order $2 \nb$. Then $\xis(t;\yb(t_0))$ is stationary and its PE decomposition is given by
\begin{align}
\label{eq:xisPE}
\xis(t;\yb(t_0)) = U \xipe(t) \,, 
\end{align}
where the component $\xipe(t) \in \RR^{n_c}$ is PE. Moreover, $\cC = \Img( R[\xis] )$ is the PE subspace of $\xis(t;\yb(t_0))$.
\tqed
\end{prop} 

Suppose that we had available to us $\xis(t;\yb(t_0))$ whose values over time span $\cC$ using a sufficiently rich input $u$. Then there are standard methods to recover $\cC$ using $\xis$.
Following ideas from \cite{UZEDA25,UZEDA23}, we consider a filter
\begin{align*}
\dot{\Sigma} &= -\veps\Sigma + \veps \xis(t;\yb(t_0)) \xis\T (t;\yb(t_0)) \,, 
\end{align*}
with $\Sigma \in \RR^{2\nb \times 2\nb}$ and $\veps > 0$.
This filter was inspired by Kreisselmeier's integral algorithm \cite{KREISSEL77} and was adopted in \cite{MARINO22} to characterize the lack of persistency of excitation of a regressor. The filter has a steady-state solution \cite[Lemma~1]{UZEDA25}, \cite[Lemma~2.1]{MARINO22} which corresponds to the PE subspace of $\xis$.

\begin{lemma} \label{lem:UZEDA}
Suppose $\xis$ is stationary and its PE decomposition is given by \eqref{eq:xisPE}. Let $\beta\pe,\, T\pe > 0$ be PE constants in \eqref{eq:PE} for $\xipe$. Consider the steady-state system
\begin{align}
\label{eq:subsest1}
\dot{\Sigma} 
&= -\veps \Sigma + \veps \xis(t;\yb(t_0)) \xis\T (t;\yb(t_0)) 
\end{align}
with initial conditions $\Sigma(t_0) = \veps T\pe R[\xis]$.
There exists a bounded symmetric $\Lambda(t) \in \RR^{n_c \x n_c}$ such that 
\begin{align*}
\Sigma(t) &= U \Lambda(t) U\T \,, \quad
\Lambda(t) \succeq \veps T\pe \beta\pe {\rm e}^{- \veps T\pe} I
\end{align*}
for all $t \geq t_0 \geq 0$. 
Hence, $\Img(\Sigma(t)) = \Img(U)$ for all $t \geq t_0 \geq 0$.
\tqed
\end{lemma}
This lemma says that the image of the filter state $\Sigma(t)$ evolves in the controllable subspace of the overparameterized model, assuming  appropriate initial conditions $\Sigma(t_0)$.
Since $\xis(t;\yb(t_0))$ is unknown, we apply a filter of the form
\begin{align}\label{eq:subsest2}
\dot{\Sigh} & = -\veps \Sigh + \veps \xih \xih^{\T} \,,
\end{align}
with $\Sigh \in \RR^{2\nb \times 2\nb}$ and $\veps > 0$. 
Define the tracking error $\Sigt := \Sigh - \Sigma(t)$, with $\Sigma(t_0)$ given in Lemma~\ref{lem:UZEDA}. One may verify by direct calculation
\begin{align}
\label{eq:Sigt}
\dot{\Sigt} &= -\veps \Sigt + \veps(\xit\xis\T(t;\yb(t_0)) 
+ \xis(t;\yb(t_0))\xit\T + \xit\xit\T)\,.
\end{align}  
From \eqref{eq:xit}, we know that $\xit(t) \rightarrow 0$ exponentially, and since $\xis(t;\yb(t_0))$ is bounded, we deduce that $\Sigt \rightarrow 0$ exponentially. In summary, we have shown that by applying \eqref{eq:subsest2} we obtain an estimate of $\Sigma(t)$, so that by Lemma~\ref{lem:UZEDA} we can asymptotically recover $\Img(U) = \cC$.  

\section{Undermodeled Plant}
\label{sec:underparam}

In this section we relax Assumption~\ref{assum:over} by considering the case when the plant order has been erroneously underestimated. 
\begin{assumption} \label{assum:under} 
The integer $\nb$ satisfies $0 < \nb < n$.
\tqed
\end{assumption}

The key observation is that the methods developed before using the overparameterized model \eqref{eq:overp2} are no longer valid. Indeed, the next result confirms that there does not exist $(\theta_1, \theta_2)$ satisfying Proposition~\ref{prop:thetas}. 

\begin{prop}
\label{prop:under}
Consider \eqref{eq:overp2} with $\nb < n$. Then $\cP = \emptyset$.
\tqed
\end{prop}
 
Since an overparameterized model is not available when $\nb < n$, we work directly with the observer \eqref{eq:observer}. 
Suppose that $u$ is stationary, so that $\xih$ is also stationary.
We define the plant order estimate to be
\begin{align}
\label{eq:nhat}
\nh &:= \rnk ( R[\xih] )  - \nb \,.
\end{align} 
Notice that if $\xih(t)$ is PE, then $R[\xih]$ has rank $2 \nb$ \cite[Prop~2.7.1]{SASTRY89}, in which case $\nh = \nb$. This result appears to be inconclusive; however, one may repeat the experiment with a larger $\nb$ until a gap $\nh < \nb$ is detected. This idea suggests a procedure that will be further developed in Section~\ref{sec:algorithm} to determine the true plant order.

To compute $R[\xih]$ we may use ideas about subspace estimation from Section~\ref{sec:subspace}. Suppose $u$ is a multisine signal, so that it is almost periodic. Then the steady-state solution $\xih_{ss}(t) \in \RR^{2 \nb}$ of \eqref{eq:observer} is also almost periodic and stationary \cite{HALE80}. Moreover, by \cite[Prop~3]{UZEDA25} it possesses a PE decomposition given by $\xih_{ss} = W \xih_{ss,pe}$, where $\Img(W) = \Img(R[\xih_{ss}])$ and $\xih_{ss,pe} \in \RR^{q\pe}$ is PE. We also know that for arbitrary initial conditions of \eqref{eq:observer}, $\xih(t) = \xih_{ss}(t) + \xit(t)$, where $\xit(t)$ denotes the transient which decays exponentially to zero. Finally, by the properties of the autocovariance (see \cite{UZEDA25}), $R[\xih] = R[\xih_{ss}]$. 

Now consider the filter
\begin{align*}
\dot{\Sigma} &= - \veps \Sigma + \veps \xih_{ss}(t) \xih_{ss}\T(t) \,,
\end{align*}
where $\veps > 0$ and $\Sigma(t) \in \RR^{2 \nb \times 2 \nb}$. 
Then by Lemma~\ref{lem:UZEDA} (with $n_c$ replaced by $q\pe$) it has a special solution $\Sigma(t)$ with initial condition $\Sigma(t_0) = \veps T\pe R[\xih_{ss}]$, where $(\beta\pe,T\pe)$ are PE constants of $\xih_{ss,pe}$. Once again we have that $\Img(\Sigma(t)) = \Img(W) = \Img(R[\xih_{ss}])$ for all $t \ge t_0 \ge 0$. 

Since we do not have access to $\Sigma(t)$, we can again use a filter of the form \eqref{eq:subsest2} (now with $\nb < n$). Define the tracking error $\Sigt := \Sigh - \Sigma(t)$, where $\Sigma(t)$ denotes the special solution. Then
\begin{align*}
\dot{\Sigt} &= -\veps \Sigt 
+ \veps( \xit \xih_{ss}\T(t) + \xih_{ss}(t)\xit\T + \xit \xit\T )\,.
\end{align*}  
Since $\xih_{ss}$ is bounded (because it is almost periodic) and $\xit(t) \to 0$ exponentially, we see that $\Sigt(t) \to 0$. As we found before in the overmodeled case, we may asymptotically
recover $\Img(W) = \Img(R[\xih])$ using \eqref{eq:subsest2}.

The remaining theoretical step to realize a procedure for the undermodeled case is to ensure that a sufficiently rich input guarantees a PE state of the observer. Using \eqref{eq:FboGbo}, we write \eqref{eq:observer} as
\begin{align*}
\dot{\xih} &= \Fbz \xih + \Gbz 
\begin{bmatrix} y \\ u \end{bmatrix} \,,
\end{align*}
where $\xih(t) \in \RR^{2\nb}$. Next, take the Laplace transform, utilizing the identity \cite[Lemma~1]{KREISSEL96}
\begin{align}
\label{eq:laplaceexpansion}
\nonumber
&(sI-\Fbz)\inv \Gbz \\
&=  \frac{1}{f(s)}
\begin{bmatrix}
1  & s & \cdots & s^{\nb-1} & 0 & 0 & \ldots & 0 \\
 0 & 0 & \ldots & 0 &1 & s & \cdots & s^{\nb-1} \end{bmatrix}\T 
 \,.
\end{align}
We obtain 
\begin{align}
\label{eq:xihreachable}
\xih(s) & = (sI- \Fbz)^{-1}\Gbz 
\begin{bmatrix} \frac{\beta(s)}{\al(s)} \\ 1 \end{bmatrix} u(s) 
 =: T(s) u(s)\,.
\end{align}
To demonstrate that $\xih$ is PE under a SR input, we cannot directly apply \cite[Theorem~6.3]{NARENDRA89}, as the controllability of $T(s)$ has not been established. By rewriting \eqref{eq:observer} as \eqref{eq:xihreachable},  $\xih \in \RR^{2\nb}$ can be regarded as outputs of the model $T(s)$, and
the standard approach to prove PE is \cite[Theorem 5.2.1]{IOANNOU12}.
 This result cannot be applied in our case as it cannot be established that $T(j\omega_1) \ldots T(j\omega_{2\nb})$ are linearly independent for all $\omega_1 \,, \omega_2 \,, \ldots \,, \omega_{2\nb} \in \RR$: the number of frequencies $2 \nb$ is not sufficient since $\deg (\alpha(s)) = n > \nb$. 
Linear independence of $T(j\omega_1) \ldots T(j\omega_{2\nb})$ for all choices of frequencies is only a sufficient condition for $\xih$ to be PE. One can still obtain a PE signal $\xih(t)$ for almost all choices of frequencies. 

\begin{theorem}\label{thm:underest}
    Consider observer \eqref{eq:observer} with $\nb < n$. Suppose the input $u(t) \in \RR$ is a multisine composed of $\nb$ randomly selected, distinct frequencies $\omega_1, \ldots, \omega_{\nb}$. Then, for almost every choice of these $\nb$ frequencies, the autocovariance of $\xih$ satisfies $\rnk(R[\xih]) = 2\nb$.
    \tqed
\end{theorem}

\section{Model Order Assignment Algorithm}
\label{sec:algorithm}

We consider a model order estimation algorithm that generates an estimate of the plant order as a function of the current value of $\nb$. Abusing notation (compare to \eqref{eq:nhat}) we define the plant order estimate to be
\begin{align}
\label{eq:nhat2}
\nh(t; \nb) &= \rnk ( \Sigh(t) ) - \nb \,.
\end{align} 
Based on our foregoing theoretical development, $\nh(t;\nb)$ converges to $\nb$ in finite time when $\nb < n$, and it converges to $n$ when $\nb \geq n$. Convergence of the rank in finite time will be achieved using a robust SVD computation as explained in Section~\ref{sec:examples}.
When $\nb = n$, the result is inconclusive, as we cannot distinguish the underparameterized and correct model order estimates. Thus, the procedure iteratively estimates $\nh$ using \eqref{eq:nhat2}, incrementing $\nb$ at each step and repeating the experiment until a gap $\nh < \nb$ is detected.  

\begin{algorithm} 
\begin{enumerate}[label=\arabic*.]
    \item 
    Initialize $\nb := 1$ and $\nh(t;0) = 0$.
    \item 
    Select the subspace estimator sensitivity rate $\veps > 0$, and the termination time $t_f$ for the experiment. Select $f(s)$ and define the observer matrices $F \in \RR^{\nb \times \nb}$ and $G \in \RR^{\nb}$ according to \eqref{eq:FG}.
    \item 
    Apply a multisine SR input $u(t)$ of order $2\nb$ to \eqref{eq:theplant}.
    \item 
    Compute $\xih(t) \in \RR^{2\nb}$ and $\Sigh(t) \in \RR^{2\nb \times 2\nb}$ for $t \in [t_0,t_f]$ using \eqref{eq:observer} and \eqref{eq:subsest2} with $\xih(t_0) = 0$, $\Sigh(t_0) = 0$.
    \item 
    Compute $\nh(t_f;\nb)$ using \eqref{eq:nhat2}.
    \item 
    If $\nh(t_f;\nb) = \nh(t_f;\nb - 1)$ then terminate the experiment, otherwise set $\nb := \nb + 1$ and return to Step~2.
\end{enumerate}
\caption{Model Order Assignment}
\end{algorithm}

\begin{remark}
Our procedure for incrementally estimating the plant order mirrors a similar procedure in the (discrete time) Ho-Kalman algorithm, where a Hankel matrix is constructed for increasing estimates of the plant order until the rank of two sequential Hankel matrices is identical. More sophisticated heuristics for termination of the algorithm could be developed based on a {\em dimension plot}, a plot of the estimate $\nh$ versus the upper bound $\nb$. Dimension plots will be illustrated in Section~\ref{sec:examples}. 
\tqed
\end{remark}

The correctness of the algorithm is guaranteed by collecting our previous rank estimates in Theorems~\ref{thm:order} and \ref{thm:underest}. Here we assume the noise free case with $t_f$ sufficiently long such that the rank computation in \eqref{eq:nhat2} is accurate. If $\nb < n$, Theorem~\ref{thm:underest} gives $\rnk(\Sigh(t_f)) = 2\nb$ which by \eqref{eq:nhat2} implies $\nh(t_f;\nb) = \nb$. Incrementing $\nb$, so long as $\nb < n$, the algorithm continues to increment $\nh$ by one on each iteration, so it is guaranteed to eventually reach $\nb = n$. At this point, by Theorem~\ref{thm:order}, $ \rnk ( \Sigh(t_f) ) = \dim(\cC) = 2 n$ so $\nh(t_f,n) = \nb = n$. However, we cannot yet conclude this is the true model order. When $\nb$ is incremented again such that $\nb = n+1$, $\Sigh(t_f)$ is no longer full rank, but rather $\rnk(\Sigh(t_f)) = \dim(\cC) = n + \nb = 2 \nb - 1$. Therefore, $\nh(t_f,n+1) = \rnk(\Sigh(t_f)) - \nb = n$. A fixed point in the iterations is reached, and the algorithm terminates. Thus, we have proved the following.

\begin{theorem}
Suppose we have a robust rank estimation procedure such that $\rnk(\Sigh(t_f)) = 2 \nb$ for $\nb <  n$, and $\rnk(\Sigh(t_f)) = n + \nb$ for $\nb \ge n$.
Then Algorithm~1 terminates with a correct estimate of the model order for \eqref{eq:theplant}.
\tqed
\end{theorem}

The key to proper operation of the algorithm is the finite time convergence of the rank of $\Sigh(t)$ to the correct value, which relies on a robust rank estimation procedure, discussed in Section~\ref{sec:examples}.

\section{Parameter Identification}
\label{sec:params}

Once the model order $n$ has been recovered using Algorithm~1, one may then proceed to identify the plant parameters using standard methods in adaptive control. Since $\nb = n$, we can fix $\delta(s) = 1$ in Proposition~\ref{prop:thetas}. Then the parameters $\theta = (\theta_1,\theta_2) \in \RR^{2 n}$ in \eqref{eq:thetas} become 
\begin{align} \label{eq:nominalparams}
\theta_1 = f - \al \,, \qquad \theta_2 = \beta \,.
\end{align}
If we can estimate $\theta$ then the plant parameters can be recovered. Let $\hat{\theta}$ denote an estimate of $\theta$, and let $\tht := \thh - \theta$ be the parameter estimation error.

Consider the model \eqref{eq:overp2}. By Proposition~\ref{prop:matching} it has an output matching solution $\xis(t;\yb(t_0))$ such that $y\0(t) = y(t) = \theta\T \xis(t;\yb(t_0))$ for all $t \ge t_0 \ge 0$. We work with the observer \eqref{eq:observer}, the tracking errors $\xit = \xih - \xis(t;\yb(t_0))$ and $\Sigt = \Sigh - \Sigma(t)$, and the prediction error
\begin{align} 
\label{eq:error}
e &:= \thh\T \xih - y(t) \\
\nonumber
  & = \thh\T \xih - \theta\T \xis(t;\yb(t_0))\\
\nonumber
  & = \tht\T \xih + \theta\T \xit \,.
\end{align}
We see that \eqref{eq:error} comprises a standard static error model, modulo an exponentially vanishing term (recall \eqref{eq:xit}), so we apply the gradient law 
\begin{align}
\label{eq:gradient}
\dot{\thh} &= - \gamma e \xih \,,
\end{align}
where $\gamma > 0$. Substitute \eqref{eq:error} into \eqref{eq:gradient} and recall the error dynamics \eqref{eq:xit} and \eqref{eq:Sigt}. The resulting error model is 
\begin{subequations}
\label{eq:errormodel}
\begin{align}
\label{eq:thtdot}
\dot{\tht} &= - \gamma \left( \xih(t) \xih\T(t) \right) \tht 
              - \gamma \left( \theta\T \xit \right) \xih(t) \\
\label{eq:Sigtdot}    
\dot{\Sigt} &=  -\veps \Sigt + \veps g(t,\xit) \\
\label{eq:xitdot}
\dot{\xit} &= \Fbz \xit  \,,
\end{align} 
\end{subequations}
where $\Fbz$ was defined in \eqref{eq:FboGbo} and 
\begin{align}
\label{eq:g}
g(t,\xit) 
&= \xis(t;\yb(t_0)) \xit\T + \xit \xis\T(t;\yb(t_0)) + \xit\xit\T \,.
\end{align} 
The proof of the following result is different from the classical result as it does not rely on arguments about the plant model \cite{KREISSEL77}, \cite[Theorem~5.2.4]{IOANNOU12}) but rather establishes PE of $\xis(t)$. This approach mimics model reference adaptive control where one leverages the PE properties of the reference model solutions \cite[Ch.~3]{SASTRY89}. 

\begin{theorem}    
\label{thm:params}
Consider the error model \eqref{eq:errormodel} with regressor $\xih(t) \in \RR^{2n}$ generated by the observer \eqref{eq:observer}. Suppose the input $u$ is stationary and SR of order $2n$. Then the equilibrium $( \tht\,, \Sigt \,, \xit ) = ( 0\,, 0\,, 0 )$ of \eqref{eq:errormodel} is exponentially stable on every closed ball of initial conditions.
\tqed
\end{theorem}

\section{Robustness and Examples}
\label{sec:examples}

Theorem~\ref{thm:params} implies that the equilibrium $(\tht,\Sigt,\xit) = (0,0,0)$ has global uniform asymptotic stability (GUAS) and local exponential stability (LES) (see \cite[Prop~1]{UZEDA23}). It is well-known that this implies the closed-loop system is robust to bounded perturbations. These stability arguments easily extend to $(\Sigt,\xit)$ for any $\nb$ to establish that  \eqref{eq:Sigtdot}-\eqref{eq:xitdot} are robust to bounded perturbations. More general robustness results when the regressor is not PE are found in \cite[Theorem~7]{UZEDA25}. It remains to discuss the robustness of the SVD computation.

The robust determination of the rank of a noisy matrix is a classical problem  \cite{GOLUB80} \cite{STEWART84}, with several recent heuristics proposed in the control literature \cite{BAUER01,CHOWDHARY11}. The main challenge is to distinguish relevant excitation from noise, by establishing a clear gap in the singular values. We adopt several strategies to deal with this problem, but if the Signal-to-Noise Ratio (SNR) is too low, the order estimation process breaks down.

Lemma~\ref{lem:UZEDA} identified a special solution $\Sigma(t)$ of \eqref{eq:subsest2}, which admits a decomposition $U \Lambda (t) U\T$, where $U$ is an orthogonal matrix and $\Lambda(t)$ is a bounded symmetric matrix such that $\Lambda(t) \succeq \veps \beta_{pe} T_{pe} e^{-\veps T_{pe}} \cdot I$. This tells us that $\Sigma(t)$ possesses an SVD of the form 
\begin{align*}
\Sigma (t) 
& = \begin{bmatrix} U_3(t) & U_4(t) \end{bmatrix} 
\begin{bmatrix}
D_3(t) & 0\\
0 & 0\\
\end{bmatrix}
\begin{bmatrix} U_3\T(t) \\ U_4\T(t) \end{bmatrix} \\
&= U \Lambda(t) U\T \,.
\end{align*}
By Lemma~\ref{lem:UZEDA} the singular values are bounded as: 
\begin{align*}
\sigma_{\min}(D_3(t)) &= \sigma_{\min}(\Lambda(t)) 
\geq  \veps \beta_{pe} T_{pe} e^{-\veps T_{pe}} \,.
\end{align*}
As for $\Sigh(t)$, we may also consider its SVD 
\begin{align*}
\Sigh(t) &= 
\begin{bmatrix} U_1(t) & U_2(t) \end{bmatrix}
\begin{bmatrix}
D_1(t) & 0\\
0 & D_2(t)\\
\end{bmatrix}
\begin{bmatrix} V_1\T(t) \\ V_2\T(t) \end{bmatrix} \,.
\end{align*}
Since we have shown $\Sigt(t) \rightarrow 0$, one can show via Weyl's Theorem \cite{UZEDA25} that 
\begin{subequations}
\label{eq:sigmatol}
\begin{align}
\liminf_{t\rightarrow \infty} {\sigma_{\min}(D_1(t))} 
&\geq \veps \beta_{pe} T_{pe} e^{-\veps T_{pe}} \\
\limsup_{t\rightarrow \infty}\sigma_{\max}(D_2(t)) &= 0 \,.
\end{align} 
 \end{subequations}
Note also in the underparameterized case, $\Sigma(t)$ is full rank, so the dimension of $D_2$ is zero. In this case, the lower bound derived for $D_1(t)$ in \eqref{eq:sigmatol} applies directly to the smallest singular value of $\hat{\Sigma}(t)$.
The bound \eqref{eq:sigmatol} suggests that we may adopt a threshold parameter $\stol$ to estimate the rank of $\Sigma(t)$ from $\Sigh(t)$. It is shown in \cite{UZEDA25} using an averaging argument that a tighter, more practical bound 
\begin{align*}
\Lambda(t) \succeq \left(\beta\pe - \delta\pe(\veps)\right) \cdot I
\end{align*}
can be formed, where $\delta\pe(\cdot)$ is a  class-$\cK$ function. This lower bound provides a significantly more accurate approximation of the excitation level of $\xis$ than the bound $\veps \beta\pe T\pe e^{-\veps T\pe}$. Thus, we define a threshold parameter $\stol = \beta\pe $ to quantify the excitation level and set $\veps$ sufficiently small in order to slow down the subspace estimator. The rank of $\Sigh(t_f)$ is estimated to be the number of singular values above $\stol$.

Next, we would like to obtain an estimate of $\stol$. The key observation is that the vector $(0,\xistwo(t))$ extracted from the output matching solution $\xis(t;\yb(t_0)) = (\xisone(t),\xistwo(t))$ lies in $\cC$, even if $\xis(t;\yb(t_0))$ has been pushed outside of $\cC$ because of noise. Therefore, a proxy for $\stol$ can be generated from the excitation in $\xistwo$ alone. Accordingly, we run a second subspace estimator driven by $\xih_2$, 
\begin{align*}
    \dot{\Sigh}_2 &= -\veps \Sigh_2 + \veps \xih_2 \xih_2^{\T} \,.
\end{align*}
An estimated, heuristic threshold is:
\begin{align*}
    \shattol \approx 0.01 \cdot \lambda_{\min}(\hat{\Sigma}_2(t_f))\,,
\end{align*}
where $t_f$ is the termination time of the experiment.
This estimate accounts for the possibility that the plant attenuates excitation in $\xih_1$ (by an unknown amount), such that for a singular value to be counted as part of the relevant excitation in $\cC$, its contribution to the signal energy must be at least $1\%$ of the excitation level of $\xih_2$. This adhoc margin was selected to be sufficiently large to noise excitation yet small enough to capture true input excitation. 

Next, we demonstrate the algorithm's performance over various relevant scenarios. To guarantee a SR input in all simulations, we follow an approach similar to \cite{HUANG22,PINTELON98} using a multisine input formed of $4\nb$ frequencies evenly spaced over two bands,  $[0.05,6]$ rad$/s$ and $[9,50]$ rad$/s$, where $\nb$ is fixed for each experiment. As per Remark~\ref{rem:FGf}, we select a nominal value of $-5$ for the roots of $f(s)$, while $t_f$ is chosen so that the system has reached steady state.

\subsection{Example with Perfect Measurements}

\begin{figure*}
  \centering
  \begin{subfigure}{0.33\linewidth}
    \includegraphics[width=\linewidth]{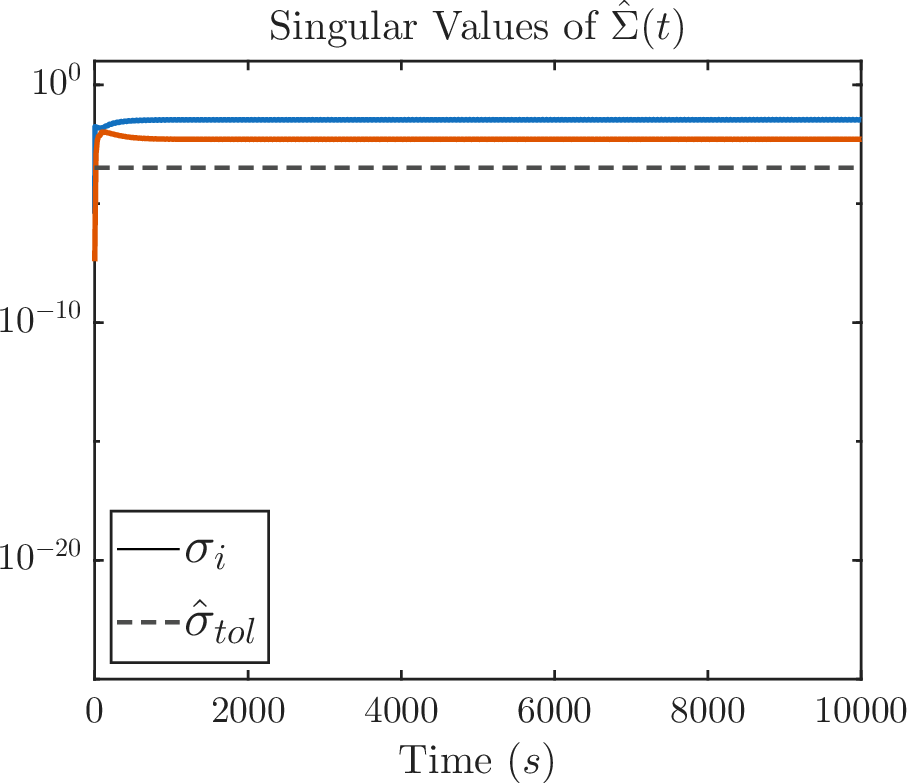}
    \caption{$\nb = 1$}\label{fig:1}
  \end{subfigure}\hfill
  \begin{subfigure}{0.33\linewidth}
    \includegraphics[width=\linewidth]{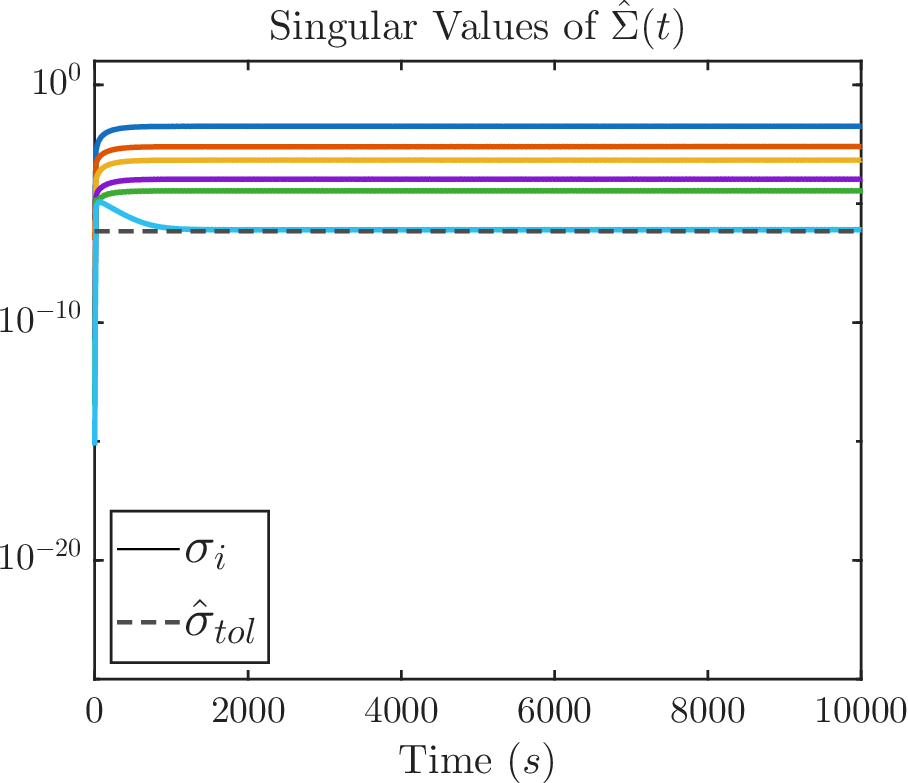}
    \caption{$\nb = 3$}\label{fig:2}
  \end{subfigure}\hfill
  \begin{subfigure}{0.33\linewidth}
    \includegraphics[width=\linewidth]{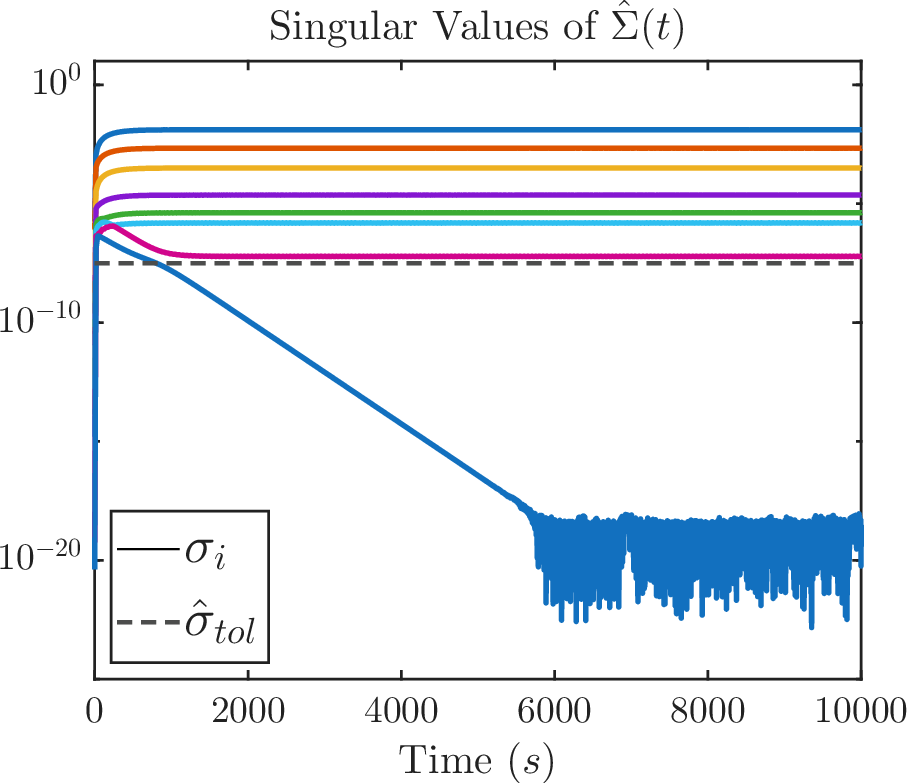}
    \caption{$\nb = 4$}\label{fig:3}
  \end{subfigure}
\caption{Evolution of the singular values of $\Sigh(t)$ (solid lines) on a log scale as a function of time, for different values of $\nb$. The value corresponding to $\shattol$, the estimate of $\stol$ (indicated by dashed lines), serves as a coarse threshold for estimating the rank of $\Sigh$.}
  \label{fig:SVDn3}
\end{figure*}

We consider a simple numerical example to illustrate the basic operation of the algorithm. Thus, consider the third-order system 
\begin{align}\label{eq:thirdorder}
    & A = \begin{bmatrix}
        -6 & -11 & -6\\
        1  &  0   & 0\\
        0  &  1   & 0
    \end{bmatrix}\,,
    B = \begin{bmatrix}
        1\\
        0\\
        0\\
    \end{bmatrix}\,,\\
   \nonumber & C = \begin{bmatrix}
       6 & 0 & 0\\
    \end{bmatrix}\,, D = 0. 
\end{align}
We first ran Algorithm~1 with $\veps = 0.005$, $t_f = 10000$s, and zero initial conditions for \eqref{eq:observer} and \eqref{eq:subsest2}. Each sinusoidal component of the input has a magnitude of one. For the sake of illustration, we computed the SVD every second, but in practice the SVD need only be computed at $t_f$, the end of each experiment. The rank estimates in Figure~\ref{fig:SVDn3} are: 
(a) $\nh = \rnk(\Sigh(t_f)) - \nb = 2 - 1 = 1$, 
(b) $\nh = \rnk(\Sigh(t_f)) - \nb = 6 - 3 = 3$, 
(c) $\nh = \rnk(\Sigh(t_f)) - \nb = 7 - 4 = 3$, at which point the algorithm terminates, since $\nb > \nh$.
We see clearly in Figure~\ref{fig:3} that one of the singular values drops to a value near machine precision. 

\begin{figure}
    \centering
    \includegraphics[width=\linewidth]{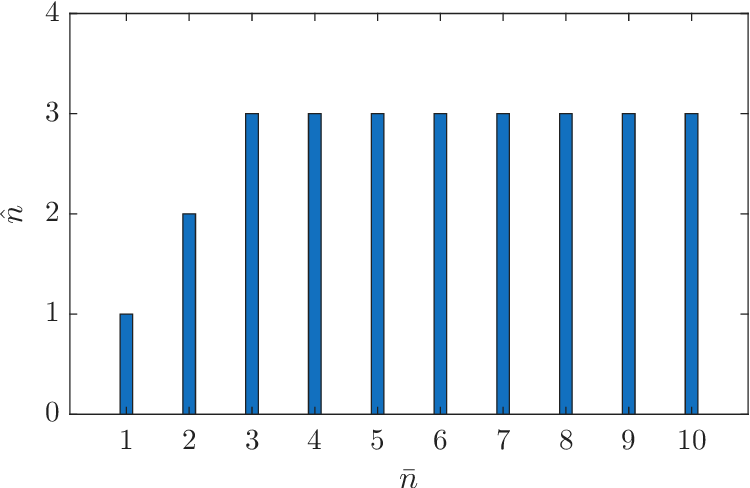}
    \caption{Dimension plot for the third-order system \eqref{eq:thirdorder}. The algorithm reaches a fixed point when $\nh(t_f;3) = \nh(t_f;4) = 3$. We generate estimates $\nh$ for larger $\nb$'s for a more complete dimension plot. The plateau in the plot provides further validation of the plant order estimate.}
    \label{fig:n3dimplot}
\end{figure}

The termination of Algorithm~1 is based on reaching a fixed point in the model order estimate. One can also generate further estimates of the plant order for larger $\nb$ to generate a {\em dimension plot}, as seen in Figure~\ref{fig:n3dimplot}. The plateau in this plot provides a visual confirmation that the true plant order has been estimated.


\subsection{Higher Order Model}

\begin{figure}[!htbp]
      \centering
    \includegraphics[width=\linewidth]{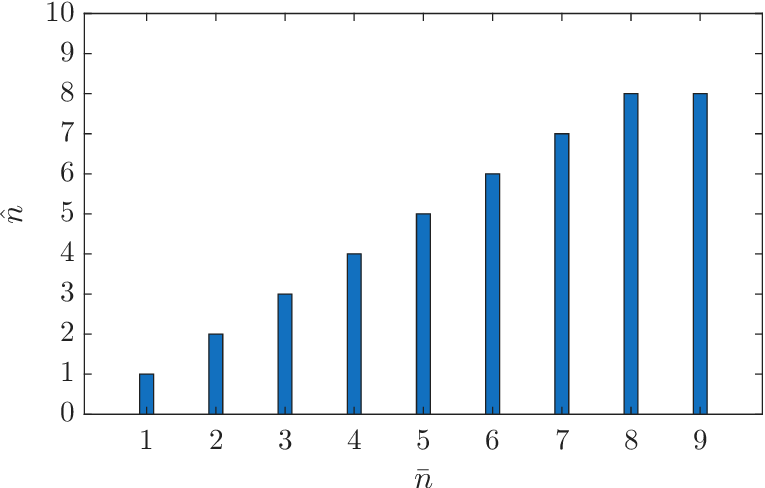}
    \caption{Dimension plot for an eighth-order system.}
    \label{fig:dimplotn8}
\end{figure}

We consider a higher-order numerical example to assess the algorithm's performance when the plant order is larger. We take an eighth-order stable SISO system with poles at
\begin{align}\label{eq:n8}
    -0.5 \pm 1j,\quad -1.0 \pm 3j,\quad -1.5 \pm 4j,\quad -2.0 \pm 6j \,.
\end{align}
 We run Algorithm~1 using $\varepsilon = 0.005$ and $t_f =15000$s. The magnitude of the input excitation gets attenuated for higher-order derivatives as $\nb$ increases. Consequently, we scale the input amplitude by a factor of $10$ to ensure that the excitation effectively propagates to the higher-order derivatives in \eqref{eq:observer} as $\nb$ increases. The resulting dimension plot is shown in
Figure~\ref{fig:dimplotn8}. The model order is thus estimated as $\nh = 8$ since we observe a plateau at $\nb = 9$.

\subsection{Robustness to Noise}

\begin{figure}[!htbp]
    \centering
    \includegraphics[width =0.9\linewidth]{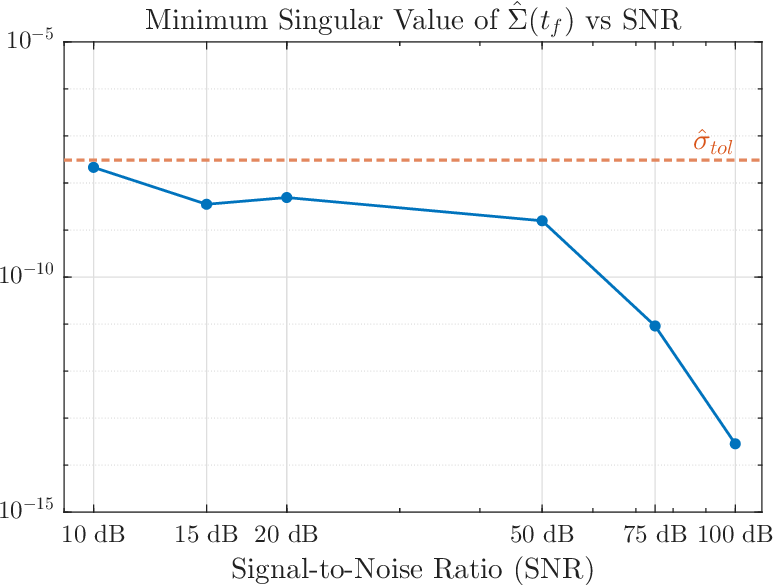}
    \caption{The mean of the minimum singular value of $\Sigh$ at steady state  plotted on a log scale. The margin below $\shattol$ decreases as the SNR increases.}
    \label{fig:ExMargin}
\end{figure}

\begin{figure}[!htbp]
  \centering

  \begin{subfigure}{0.9\linewidth}
    \centering
    \includegraphics[width=\linewidth]{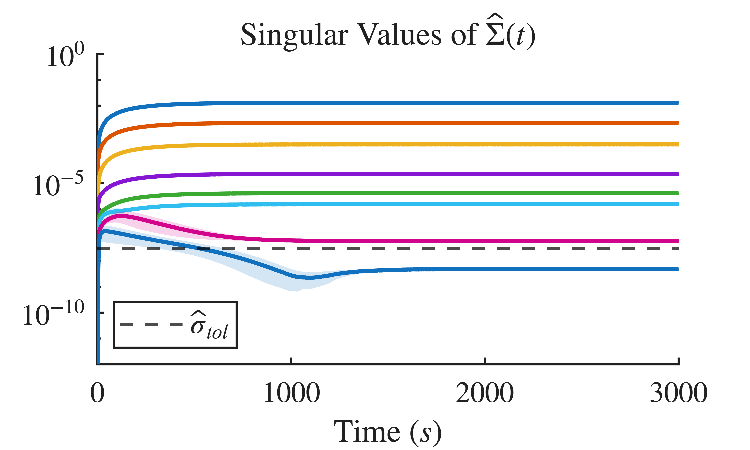}
    \caption{SNR $ = 20 \mathrm{dB}$}
    \label{fig:SNR20}
  \end{subfigure}

  \vspace{0.6em}

  \begin{subfigure}{0.9\linewidth}
    \centering
    \includegraphics[width=\linewidth]{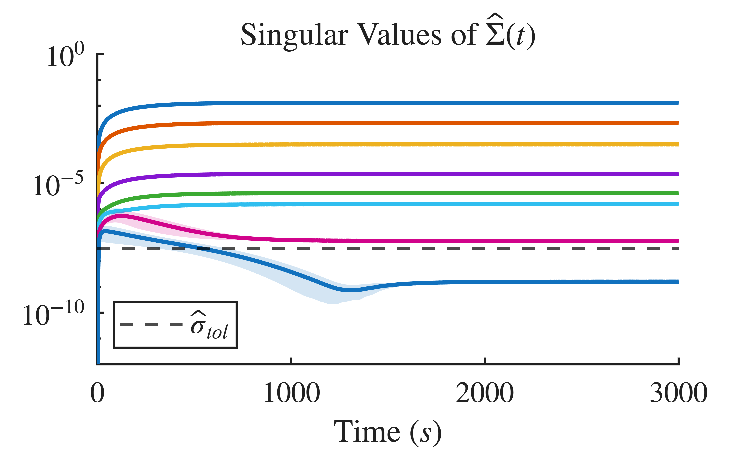}
    \caption{SNR $= 50 \mathrm{dB}$}
    \label{fig:SNR50}
  \end{subfigure}

  \caption{Evolution of the mean of singular values of $\Sigh(t)$ (solid lines) in the presence of measurement noise for (a) SNR$=20\mathrm{dB}$ and (b) SNR$=50\mathrm{dB}$, plotted on a log scale. The shaded regions correspond to $\pm 1$ standard deviation bands across $50$ experiments. The gap in excitation below $\shattol$ increases with increasing SNR from $20\mathrm{dB}$ to $50\mathrm{dB}$.}
  \label{fig:NoiseOrders}
\end{figure}

Consider again the system \eqref{eq:thirdorder} subject to zero-mean measurement noise at different SNRs at the plant output. We consider an experiment with $\nb=4$, $\veps=0.005$, and $t_f = 5000$s. Since the true system order is $n = 3$, $\sigma_{\max}(D_2(t))$ in \eqref{eq:sigmatol} corresponds to the minimum singular value of $\Sigh(t)$. As seen in \eqref{eq:sigmatol}, this singular value vanishes in the noise-free case. 
Now we are interested to investigate the effect of noise on this minimum singular value. Figure~\ref{fig:ExMargin} shows the mean of the minimum singular value of $\Sigh(t_f)$, across $50$ Monte Carlo trajectories, for each choice of SNR. The estimate $\shattol$ is the same over all experiments, implying it does not pick up noise. On the other hand, we see that as the SNR at the plant output decreases, the gap in excitation level above noise becomes smaller, making the model order increasingly difficult to identify. 

Figure~\ref{fig:NoiseOrders} provides a different perspective of the same experiment, by showing the average (solid line) and variance bands (shaded area) of each singular value as a function of time, with SNR of $20\mathrm{dB}$ and $50 \mathrm{dB}$. We see in Figure~\ref{fig:NoiseOrders}(a) that the excitation gap below the threshold is smaller at $20\mathrm{dB}$ than at $50\mathrm{dB}$ in (b). 
Nevertheless, we observe that the variance of the singular values of $\Sigh(t)$ is small such that the rank estimate is consistent across all Monte Carlo runs for each SNR.

\section{Proofs}

\subsection{Proof of Lemma~\ref{lem:PEC}}

First suppose $\cC(A,B) = \{ 0 \}$, implying $B = 0$. For any initial condition, $x(t)$ is stationary, and $x(t) \rightarrow 0$, which implies $R[x] = 0$. Hence, $\Img(R[x]) = \{ 0 \} = \cC(A,B)$.
Second, suppose $\cC(A,B) = \RR^n$. Since $u$ is stationary and SR of order $n$, and $A$ is Hurwitz, by \cite[Theorem~6.3]{NARENDRA89}, $x(t)$ is PE, and its autocovariance $R[x]$ is positive definite by \cite[Prop~2.7.1]{SASTRY89}. 
Hence, $\Img{(R[x])} = \RR^n = \cC(A,B)$.

Finally, suppose $0 < n_1 := \dim (\cC(A,B)) < n$. 
Let $\begin{bmatrix} U & U\pr \end{bmatrix} \in \RR^{n \times n}$ be an orthogonal matrix such that $\Img(U) = \cC(A,B)$ and $\Img(U\pr) = \cC(A,B)^{\perp}$. Define the coordinate transformation
\begin{align}
\label{eq:Kalman0}
\begin{bmatrix} x_c \\ x_u \end{bmatrix} 
&= \begin{bmatrix} U & U\pr \end{bmatrix}\inv x 
 = \begin{bmatrix} U\T \\ U\pr\T \end{bmatrix} x \,.
\end{align}
According to the Kalman decomposition, 
\begin{subequations}
\label{eq:Kalman1}
\begin{align}
\label{eq:xc}
\dot{x}_{c} &= A_{11} x_c + A_{12} x_u + B_1 u \\
\label{eq:xu}
\dot{x}_{u} &= A_{22} x_u \,,
\end{align}
\end{subequations}
where $A_{11}$, $A_{22}$ are Hurwitz by assumption, and $x_c \in \RR^{n_1}$ and $x_u \in \RR^{n-n_1}$ denote the states of the controllable and uncontrollable subsystems, respectively. Since \eqref{eq:xc} is an LTI system, it can be split as
\begin{subequations}
\begin{align}
\label{eq:xc1}
\dot{x}_{c1} &= A_{11} x_{c1} + B_1 u\\
\label{eq:xc2}
\dot{x}_{c2} &= A_{11} x_{c2} + A_{12} x_u 
\end{align}
\end{subequations}
such that $x_c(t) = x_{c1}(t) + x_{c2}(t)$ under a suitable choice of initial conditions $(x_{c1}(0), x_{c2}(0) )$. Since $A_{11}$ and $A_{22}$ are Hurwitz, $x_u(t) \rightarrow 0$ and $x_{c2}(t) \rightarrow 0$ exponentially. Since the pair $(A_{11}, B_1)$ is controllable by the Kalman decomposition, $A_{11}$ Hurwitz, and $u$ is stationary and SR of order $n \ge n_1$, we can apply \cite[Theorem~6.3]{NARENDRA89} to obtain $x_{c1}$ is PE. Thus, 
$x_c(t) = x_{c1}(t) + x_{c2}(t)$ with $x_{c2}(t) \rightarrow 0$, so by \cite[Lemma~6.5(ii)]{NARENDRA89}, $x_c$ is PE.

By assumption input $u$ is stationary and SR of order $n$, so by \cite[Prop~1.6.2]{SASTRY89}, we know $x(t)$ and $x_c(t)$ are stationary. Since $x_c$ is also PE, 
by \cite[Prop~2.7.1]{SASTRY89}, $R[x_c]$ is positive definite. Observe from \eqref{eq:Kalman0} that 
\begin{align*}
x = U x_c + U\pr x_u \,.
\end{align*} 
Using the fact that $x_u(t) \to 0$, we have
$R[x] = R[ U x_c] = U R[x_c] U\T$.
Because $R[x_c]$ is positive definite, it is straightforward to show
$\Img( U R[x_c] U\T ) = \Img(U)$. We conclude that 
\begin{align*}
\cC(A,B) &= \Img(U) = \Img(R[x]) \,.
\end{align*}

\subsection{Proof of Proposition~\ref{prop:thetas}}

First suppose $\theta$ satisfies \eqref{eq:thetas}. We will show $\theta \in \cP$. Notice that because of \eqref{eq:FG}, 
\begin{align}
\label{eq:derivatives}
(sI - F)\inv G &= 
\frac{1}{f(s)} 
\begin{bmatrix} 1 & s & \cdots & s^{\nb-1} \end{bmatrix}\T \,.
\end{align}
Taking Laplace transforms, \eqref{eq:overp1} becomes 
\begin{align*}
\xi_1(s) &= (sI - F)\inv G y\0(s) \\
\xi_2(s) &= (sI - F)\inv G u(s)   \\
y\0(s)   &= \theta_1\T \xi_1(s) + \theta_2\T \xi_2(s) \\
         &= \frac{1}{f(s)} (f(s) - \al_e(s)) y\0(s) + \frac{1}{f(s)} \beta_e(s) u(s) \,.
\end{align*}
Multiplying by $f(s)$ on both sides and canceling terms, we obtain 
\begin{align*}
\frac{y\0(s)}{u(s)} &= \frac{\beta_e(s)}{\alpha_e(s)} 
                     = \frac{\beta(s)}{\alpha(s)} \,.
\end{align*}
Thus, $\theta \in \cP$.

For the converse direction, suppose $\theta \in \cP$. Then 
\begin{align}
\label{eq:match1}
\theta\T (sI - \Fb(\theta))\inv \Gb = \frac{\beta(s)}{\alpha(s)} \,.
\end{align}
Using \eqref{eq:overp2} and \eqref{eq:match1}, we have
\begin{align}
\nonumber
(sI - \Fb(\theta))\inv \Gb 
&=
\begin{bmatrix}
(sI-F)\inv G \theta\T (sI-\Fb(\theta))^{-1}\Gb  \\
(sI-F)\inv G
\end{bmatrix} \\
\label{eq:match2}
&= 
\begin{bmatrix}
(sI-F)\inv G \frac{\beta(s)}{\alpha(s)} \\
(sI-F)\inv G
\end{bmatrix} \,.
\end{align}
Applying \eqref{eq:derivatives} and \eqref{eq:match2} to \eqref{eq:match1}, we obtain
\begin{align*}
\frac{\beta(s)}{\alpha(s)}
&= \theta_1\T (sI-F)\inv G \frac{\beta(s)}{\alpha(s)} + 
   \theta_2\T (sI-F)\inv G \\
&= 
\frac{\theta_1(s)}{f(s)} \frac{\beta(s)}{\alpha(s)} + \frac{ \theta_2(s)}{f(s)} \,.
\end{align*}
Rearranging, this becomes
\begin{align*}
(f(s) - \theta_1(s) ) \beta(s) &= \theta_2(s) \alpha(s) \,.
\end{align*}
Now $\deg(f(s)) = \nb$ and $f(s)$ is monic, and $\deg(\theta_1(s)) \le \nb -1$. Therefore, $f(s) - \theta_1(s)$ has degree $\nb$ and it is monic (for all choices of $\theta_1(s)$). Also, $\al(s)$ is monic with degree $n$.
Since $\{ \beta(s), \alpha(s) \}$ are coprime, then $\al(s)$ is a factor of $f(s) - \theta_1(s)$. This means there must exist a monic polynomial $\delta(s)$  of degree $p = \nb - n$ such that 
\begin{align*}
f(s) - \theta_1(s) &= \delta(s) \alpha(s) \,, \quad
\theta_2(s)         = \delta(s) \beta(s) \,.
\end{align*}
This is the same as \eqref{eq:thetas}, which proves the result.

\subsection{Proof of Proposition~\ref{prop:matching}}
\label{proof:prop:matching}

Fix $t_0 \geq 0$.
The output response of \eqref{eq:overp2}  is:
\begin{align*}
y\0(t) &=  \theta\T e^{\Fb(\theta) (t-t_0)}\xi(t_0)+ \int^t_{t_0} 
\theta\T e^{\Fb(\theta) (t-\tau)}\Gb u(\tau) d\tau 
\end{align*}
for all $t\geq t_0 \ge 0$. From Proposition~\ref{prop:thetas} we know
\begin{align*}
\theta\T (sI - \Fb(\theta) )\inv \Gb 
&= C (sI- A)\inv B = \frac{\beta(s)}{\al(s)} \,,
\end{align*}
where $(C,A,B)$ is a minimal realization of \eqref{eq:theplant}.
Taking the inverse Laplace transform yields
\begin{align}
\label{eq:impulse}
\theta\T {\rm e}^{\Fb(\theta) t} \Gb &= C {\rm e}^{At} B \,, 
\qquad \forall t \ge 0 \,. 
\end{align}
Hence, we have for all $t \ge t_0 \ge 0$,
\begin{align}\label{eq:yoProp1}
y\0(t) 
&=  \theta\T e^{\Fb(\theta) (t-t_0)}\xi(t_0) + \int^t_{t_0} C {\rm e}^{A(t-\tau)} B u(\tau) d\tau \,.
\end{align}

Next we address the initial conditions. Taking repeated derivatives of \eqref{eq:impulse} we find 
\begin{align}
\label{eq:impulse2}
\theta\T {\rm e}^{\Fb(\theta) t} \Fb(\theta)^i \Gb 
&= C {\rm e}^{A t} A^i B \,,
 \qquad i \ge 0 \,, \quad \forall t \ge 0 \,.
\end{align}
Let $x(t_0)$ be the initial condition of \eqref{eq:ABC}. Since $(A,B)$ is controllable, we can find $c_0 \in \RR^n$ such that
\begin{align*}
x(t_0) &= 
\begin{bmatrix} B & A B & \cdots & A^{n-1} B \end{bmatrix} c_0 \,.
\end{align*}
Then select the initial condition 
\begin{align}
\label{eq:ic}
\xi(t_0) &:= 
\begin{bmatrix} \Gb & \Fb(\theta) \Gb & \cdots & \Fb(\theta)^{n-1} \Gb \end{bmatrix} 
c_0 \,.
\end{align}
Notice that $\xi(t_0) \in \cC$, implying $\xi(t) \in \cC$ for all $t \ge t_0 \ge 0$, which proves (ii). To prove (i), we use \eqref{eq:impulse2} to obtain
\begin{align*}
\theta\T &{\rm e}^{\Fb(\theta) (t - t_0)} \xi(t_0) \\
&= 
\theta\T {\rm e}^{\Fb(\theta) (t - t_0)}
\begin{bmatrix} \Gb & \Fb(\theta) \Gb & \cdots & \Fb(\theta)^{n-1} \Gb \end{bmatrix} 
c_0 \\
&= 
C {\rm e}^{A (t - t_0)}
\begin{bmatrix} B & A B & \cdots & A^{n-1} B \end{bmatrix} 
c_0 \\
&= C {\rm e}^{A (t - t_0)} x(t_0) \,.
\end{align*}
Therefore, \eqref{eq:yoProp1} becomes
\begin{align*}
y\0(t) &=  C e^{A(t-t_0)}x(t_0)
+ \int^t_{t_0} C {\rm e}^{A(t-\tau)} B u(\tau) d\tau
= y(t) \,, 
\end{align*}
for all $t \ge t_0 \ge 0$. 

\subsection{Proof of Proposition~\ref{prop:controllable} }

(i) Let $\theta, \theta' \in \cP$. By definition of $\cP$
\begin{align}
\label{eq:match3}
(\theta')\T (sI-\Fb(\theta'))\inv\Gb &= \theta\T (sI-\Fb(\theta))\inv \Gb
                               = \frac{\beta(s)}{\alpha(s)} \,. 
\end{align}
Also from \eqref{eq:match2}, $\theta', \theta \in \cP$ implies
\begin{align}
\label{eq:match4}
(sI - \Fb(\theta'))\inv \Gb &= (sI - \Fb(\theta))\inv \Gb \,.
\end{align}
Taking the inverse Laplace transform on both sides,
\begin{align*}
{\rm e}^{\Fb(\theta') t} \Gb &= {\rm e}^{\Fb(\theta) t} \Gb \,,  
\quad  \forall t \ge 0 \,.
\end{align*}
Next we take repeated derivatives and evaluate at $t = 0$ to get
\begin{align}
\label{eq:Cmatch}
\left( \Fb(\theta') \right)^i \Gb &= \left( \Fb(\theta) \right)^i \Gb \,, 
\quad i = 0, 1, \ldots 
\end{align}
Hence, $\cC(\theta') = \cC(\theta)$, as desired.

(ii) Let $\theta, \theta' \in \cP$ and select $\xi(t_0;\theta)$ and $\xi(t_0;\theta')$ according to Proposition~\ref{prop:matching}. 
Using \eqref{eq:ic} and \eqref{eq:Cmatch}, we have that
\begin{align*}
\xi(t_0;\theta) 
&= \begin{bmatrix} \Gb & \Fb(\theta) \Gb & \cdots 
                        & \Fb(\theta)^{n-1} \Gb \end{bmatrix} c_0 \\
&= \begin{bmatrix} \Gb & \Fb(\theta') \Gb & \cdots 
                        & \Fb(\theta')^{n-1} \Gb \end{bmatrix} c_0 \\
&= \xi(t_0;\theta') \,.
\end{align*}
Also by Proposition~\ref{prop:matching}, $y(t) = \theta\T \xis(t;\theta,\yb(t_0)) = (\theta')\T \xis(t;\theta',\yb(t_0))$ for $t \ge t_0 \ge 0$. Thus, $\xis(t;\theta,\yb(t_0))$ and $\xis(t;\theta',\yb(t_0))$ have the same initial condition and they evolve according to the same dynamics \eqref{eq:overp1}, so they are identical.

\subsection{Proof of Theorem~\ref{thm:order}}

Based on Proposition~\ref{prop:controllable} we are free to select  $\delta(s)$ to allow our analysis of the plant order to go forward with greatest ease. As such, we assume $\delta(s)$ is any monic polynomial of degree $p$ whose roots are distinct and 
\begin{align}
\label{as:delta}
(\roots(\alpha(s)) \cup \eig(F)) \cap \roots(\delta(s)) &= \emptyset \,.
\end{align}

To ease the notation we write $\Fb \equiv \Fb(\theta)$.
Define the block matrices
\begin{align*}
A_1 &= F + G \theta_1\T \,, \qquad 
C_1 = \theta_2\T \,.
\end{align*}
We employ the PBH test to analyze the dimension of the controllable subspace of \eqref{eq:overp2}. Consider
\begin{align}
\label{eq:PBH1}
   \begin{bmatrix}
        \Fb - \lambda I & \Gb
   \end{bmatrix} 
   &= \begin{bmatrix}
        A_1-\lambda I & G C_1 & 0\\
        0 & F-\lambda I & G
    \end{bmatrix} \,.
\end{align}
Due to the upper block diagonal structure of $\Fb$, we see that $\eig(\Fb) = \eig(A_1) \cup \eig(F)$. Recalling \eqref{eq:extended} and using the assumption on $\delta(s)$, we have 
\begin{align}
\label{eq:sigmaAe}
\eig(A_1) &= \roots(\al(s)) \sqcup \roots(\delta(s)) \,.
\end{align} 
We consider three cases:

(i) \ $\lambda \notin \eig(A_1)$: \
We have $\rnk (A_1 -\lambda I) = \nb$ since $\lambda \notin \eig(A_1)$. Also, $\rnk (\begin{bmatrix} F-\lambda I & G \end{bmatrix}) = \nb$ since $(F,G)$ is controllable. By the block upper diagonal form of $\begin{bmatrix} \Fb - \lambda I  & \Gb \end{bmatrix}$ it follows that $\rnk \begin{bmatrix} \Fb - \lambda I  & \Gb \end{bmatrix} = 2\nb$.

(ii) \ $\lambda \in \roots(\alpha(s))$: \
Consider the block form 
\begin{align}
\label{eq:block1}
    \begin{bmatrix}
    \Fb - \lambda I & \Gb
\end{bmatrix} = \begin{bmatrix}
    M_{11} & 0_{(\nb-1)\times \nb} & 0_{(\nb-1)\times 1}\\
    * & M_{22} & 0_{\nb\times 1}\\
    * & * & 1
\end{bmatrix}\,,
\end{align}    
where 
\begin{align}\label{eq:block2}
    M_{11} = \begin{bmatrix}
        -\lambda & 1 &\\
        &\ddots & \ddots\\
        &&-\lambda & 1
    \end{bmatrix} \in \mathbb{R}^{(\nb-1)\times \nb}
\end{align} and
\begin{align}\label{eq:block3}
   M_{22} =  \begin{bmatrix}
        \bar{\beta}_0 & \ldots & &&\bar{\beta}_{\nb-2} & \bar{\beta}_{\nb-1}\\
        -\lambda & 1 & 0 & \\
        0 & -\lambda & 1 & 0\\
        \vdots & & \ddots & \ddots & \ddots\\
        \\
        &&&-\lambda &1&0\\
        0 & \ldots &  &0& - \lambda & 1
    \end{bmatrix} \in \mathbb{R}^{\nb\times \nb}\,.
\end{align}
Notice that due to its structure, $M_{11}$ has rank $(\nb-1)$.
Also notice that by the special form of $M_{22}$, it has a Laplace expansion given by:
\begin{equation}\label{eq:laplace_exp}
    \det(M_{22}) = \bar{\beta}_{\nb-1} \lambda^{\nb-1} +  \bar{\beta}_{\nb-2} \lambda^{\nb-2} + \ldots + \bar{\beta}_0 = \beta_e(\lambda)\,.
\end{equation}
Details of this expansion are given in Section~\ref{sec:Laplace}. By \eqref{eq:sigmaAe} we know $\al(\lambda) = 0$ and $\delta(\lambda) \neq 0$. Suppose by way of contradiction that $\beta_e(\lambda) = \beta(\lambda) \delta(\lambda) = 0$. Since $\delta(\lambda) \neq 0$, then $\beta(\lambda) = 0$. This contradicts that $\{\al(s), \beta(s)\}$ are coprime polynomials. Hence, $\det(M_{22}) \neq 0$. 
Due to the block lower triangular form of \eqref{eq:block1}, it follows that $\rnk \begin{bmatrix} \Fb - \lambda I  & \Gb \end{bmatrix} =2\nb$.

(iii) $\lambda \in \roots(\delta(s))$: \ 
Recall the block form \eqref{eq:block1} composed of \eqref{eq:block2} and \eqref{eq:block3}, and consider \eqref{eq:laplace_exp}. If $\delta(\lambda)= 0$ then $\beta_e(\lambda) = \beta(\lambda)\delta(\lambda) = 0 $. Then $M_{22}$ is singular, so $\lambda \in \roots(\delta(s))$ is uncontrollable by the PBH test. From Assumption~\ref{as:delta}, $\delta(s)$ is a polynomial of degree $p = \nb - n$ whose roots are distinct. Thus there exist at least $\nb - n$ distinct uncontrollable eigenvalues of $\Fb$. 

Gathering all information, the first two cases tell us that for all $\lambda \in \roots(\al(s)) \cup \eig(F)$, $\lambda$ is controllable. This means at least $\nb + n$ eigenvalues of $\Fb$ are controllable. The third case tells us $\Fb$ has at least $p = \nb - n$ distinct uncontrollable eigenvalues corresponding to the roots of $\delta(s)$. Moreover, if $\lambda \in \eig(\Fb)$ is uncontrollable, then $\lambda \in \roots(\delta(s))$; otherwise by Assumption~\ref{as:delta}, $\lambda \in \roots(\al(s)) \cup \eig(F)$ is controllable. We conclude 
\begin{align*}
\dim(\cC) &= 2\bar{n} - \text{(\#uncontrollable eigenvalues of $\Fb$)} \\
         &= 2\bar{n} - p = n + \nb \,.
\end{align*}

\subsection{Proof of Corollary~\ref{cor:poles}}

If $\nb = n$, then $\cC = \RR^{2n}$, so all modes are controllable.
Suppose instead $\nb > n$. Let $\theta \in \cP$ correspond to any monic polynomial $\delta(s)$ of degree $p$ with distinct roots and satisfying \eqref{as:delta}.
By Theorem~\ref{thm:order}, the roots of $\alpha(s)$ and $f(s)$ are controllable, while the roots of $\delta(s)$ are uncontrollable.

Next consider any $\theta' \in \cP$ with corresponding $\delta'(s)$.
Based on Proposition~\ref{prop:controllable}(i), a controllable decomposition can be obtained for both $(\Fb(\theta), \Gb)$ and $(\Fb(\theta'), \Gb)$ using the same transformation 
$T = \begin{bmatrix} U & U\pr \end{bmatrix}$,  where
$\cC = \Img(U)$ and $\cCp = \Img(U\pr)$. 
As in \cite[p.~130]{KAILATH80}, one may verify
\begin{align*}
T^{-1} (sI-\Fb(\theta))^{-1} \Gb & = 
\begin{bmatrix} 
(sI-\Fb_c(\theta))^{-1}\Gb_c \\ 0
\end{bmatrix} \\
T^{-1} (sI-\Fb(\theta'))^{-1} \Gb & = 
\begin{bmatrix}
 (sI-\Fb_c(\theta'))^{-1}\Gb_c \\  0
 \end{bmatrix} \,,
\end{align*}
where $(\Fb_c(\theta),\Gb_c)$ and $(\Fb_c(\theta'),\Gb_c)$  are the controllable subsystems with $\Fb_c(\theta)$, $\Fb_c(\theta') \in \RR^{n_c \times n_c}$, where $n_c = \dim(\cC)$. Using \eqref{eq:match4}, this implies 
\begin{align*}
    (sI-\Fb_c(\theta))^{-1}\Gb_c = (sI-\Fb_c(\theta'))^{-1}\Gb_c\,.
\end{align*}
Taking the inverse Laplace transform, we have
\begin{align*}
{\rm e}^{\Fb_c(\theta)t}\Gb_c 
&= {\rm e}^{\Fb_c(\theta')t}\Gb_c\,, \quad t \geq 0 \,.
\end{align*}
Next, we take repeated derivatives and evaluate at $t = 0$ to obtain
\begin{align}\label{eq:MarkovEqual}
\Fb_c (\theta)^i \Gb_c =  \Fb_c (\theta') ^i\Gb_c \,, \quad i = 0\,,1\,, \ldots
\end{align}
Define the $n_c \times n_c$ controllability matrix
\begin{align*}
    Q_c :&= \begin{bmatrix}
    \Gb_c & \Fb_c(\theta) \Gb_c & \ldots &  \Fb_c(\theta)^{n_c-1} \Gb_c
    \end{bmatrix} \,.
\end{align*}
Left multiplying $Q_c$ by $\Fb_c(\theta)$ and using \eqref{eq:MarkovEqual}
we have
\begin{align*}
\Fb_c(\theta) Q_c &= \Fb_c(\theta') Q_c \,.
\end{align*}
Since $Q_c$ is invertible, $\Fb_c(\theta) = \Fb_c(\theta')$. 
We conclude, $(\Fb_c(\theta), \Gb_c) = (\Fb_c(\theta'), \Gb_c)$, so they have the same controllable eigenvalues.  
Hence, the controllable modes of \eqref{eq:overp2} using $\theta' \in \cP$ are the roots of $\alpha(s)$ and $f(s)$. Any additional modes, namely the roots of $\delta'(s)$, are uncontrollable.

\subsection{Proof of Proposition~\ref{prop:PE}}

By Theorem~\ref{thm:order} we know $n_c = \nb + n$. First suppose $\nb = n$, so that $n_c = 2 n$. For $\theta \in \cP$ (where $\cP$ is now a singleton), we have $(\Fb(\theta),\Gb)$ is controllable. Moreover,  $\theta_1 = f - \alpha$, $\theta_2 = \beta$, and 
\begin{align*}
\Fb(\theta) &=
\begin{bmatrix}
(A\0 - G \alpha\T) & G \beta\T \\
0                  & F
\end{bmatrix}
\end{align*}
is a Hurwitz matrix. Assuming $u$ is SR of order $2 n$, any solution of \eqref{eq:overp2} is PE, according to \cite[Theorem~6.3]{NARENDRA89}. In particular, $\xis(t;\yb(t_0))$ is PE.

Second, suppose $\nb > n$, so that $0 < n_c < 2 \nb$.
Apply the coordinate transformation 
\begin{align}
\label{eq:transform}
\begin{bmatrix} \xi_c \\ \xi_u \end{bmatrix}
&= \begin{bmatrix} U & U\pr \end{bmatrix}\inv \xi 
 = \begin{bmatrix} U\T \\ U\pr\T \end{bmatrix} \xi \,,
\end{align}
where $\xi_c \in \RR^{n_c}$ is the state of the controllable subsystem
and $\xi_u \in \RR^{2 \nb - n_c}$ is the state of the uncontrollable subsystem. Then \eqref{eq:overp2} becomes
\begin{subequations}
\label{eq:Kalman}
\begin{align}
\label{eq:xic}
\dot{\xi}_{c} &= F_{11} \xi_c + F_{12} \xi_u + G_1 u \\
\label{eq:xiu}
\dot{\xi}_{u} &= F_{22} \xi_u \,.
\end{align}
\end{subequations}
Since $\xis$ lies in $\cC$, we have 
\begin{align*}
\xis(t;\yb(t_0)) &= U \xi_c(t) \,, \qquad t \ge t_0 \ge 0 
\end{align*}
and $\xi_c$ evolves according to
\begin{align*}
\dot{\xi}_{c} &= F_{11} \xi_c + G_1 u \,.
\end{align*}
By the Kalman decomposition, $(F_{11},G_1)$ is a controllable pair.
By Corollary~\ref{cor:poles}, the controllable modes of \eqref{eq:overp2} are stable, so $F_{11}$ is Hurwitz. Since $u$ is SR of order $2 \nb > n_c$, we can apply \cite[Theorem~6.3]{NARENDRA89} to obtain $\xi_{c}$ is PE. 
By assumption, input $u$ is stationary and SR of order $2 \nb$. Also, 
$\xis(t;\yb(t_0))$ and $\xi_c(t)$ and bounded and piecewise continuity. 
Then by \cite[Prop~1.6.2]{SASTRY89},
$\xis(t;\yb(t_0))$ and $\xi_c(t)$ are stationary. 
In particular $R[\xis]$ and $R[\xi_c]$ exist. By definition, $\Img ( R[\xis] )$ is the PE subspace of $\xis$, and by Proposition~\ref{prop:PEdecomp}, its PE decomposition is given by \eqref{eq:xisPE}.
Since $\xi_c$ is PE, by \cite[Prop~2.7.1]{SASTRY89}, $R[\xi_c]$ is positive definite. Observe that $R[\xis] = R[ U \xi_c] = U R[\xi_c] U\T$. 
Because $R[\xi_c]$ is positive definite, it is straightforward to show
$\Img( U R[\xi_c] U\T ) = \Img(U)$. We conclude that 
$\cC = \Img(U) = \Img(R[\xis])$. Finally, we associate $\xi_c(t)$ with $\xi\pe(t)$ to obtain \eqref{eq:xisPE}.

\subsection{Proof of Proposition~\ref{prop:under}}
We know that $\deg( \beta(s)) \le n-1$, $\deg( \al(s) ) = n$, and $\deg(f(s)) = \nb$. Consider \eqref{eq:overp1} and let 
\begin{align*}
    \theta_1 := (\theta_{1,0},\ldots, \theta_{1,\nb - 1}), \qquad
    \theta_2 := (\theta_{2,0},\ldots, \theta_{2,\nb - 1})\,.
\end{align*} 
Retracing Proposition~\ref{prop:thetas}, we reach
\begin{align*}
   \frac{y\0(s)}{u(s)} & = \left [ (\theta_1)\T \frac{\beta(s)}{\al(s) f(s)} + (\theta_2)\T \frac{1}{f(s)} \right] 
   \begin{bmatrix} 1 & s & \cdots & s^{\nb-1} \end{bmatrix}\T \\
&= \frac{ \theta_1(s) \beta(s) }{ \al(s) f(s) } + \frac{\theta_2(s)}{f(s)} \,,
\end{align*}
where
\begin{align*}
\theta_1(s) &= \theta_{1,\nb - 1}s^{\nb - 1} + \theta_{1,\nb - 2}s^{\nb - 2} + \ldots + \theta_{1,0}\,,  \\
\theta_2(s) &= \theta_{2,\nb - 1}s^{\nb - 1} + \theta_{2,\nb - 2}s^{\nb - 2} + \ldots + \theta_{2,0}\,.
\end{align*}
Notice that $\deg(\theta_1(s)), \deg(\theta_2(s)) \leq \nb - 1$. Suppose by way of contradiction that 
\begin{align*}
\frac{\theta_1(s)\beta(s) + \theta_2(s) \alpha(s)}{\alpha(s)f(s)} 
&= \frac{\beta(s)}{\alpha(s)} \,.
\end{align*}
Simplifying, we have
$\theta_2(s) \alpha(s) = \beta(s) ( f(s) - \theta_1(s) )$.
Because $\{ \alpha(s), \beta(s) \}$ are coprime, the latter equality
can only hold if every root of $\al(s)$ is a root of $( f(s) - \theta_1(s) )$. However, $\al(s)$ has $n$ roots, whereas $( f(s) - \theta_1(s) )$ has $m \le \nb < n$ roots, a contradiction. 

\subsection{Proof of Theorem~\ref{thm:underest}}

Consider \eqref{eq:observer} with $\xih \in \RR^{2\nb}$ and $1 \le \nb < n$. Suppose $u$ is a multisine input formed of distinct, randomly selected frequencies $\omega_i$, $i = 1\ldots \nb$, such that $u$ is stationary and SR of order $2\nb$. We want to show that for almost all such $u$, $\xih$ is PE.
Based on \cite[Prop~2.7.1]{SASTRY89}, we show that $R[\xih]$ is positive definite.
The spectral measure $S_u(\omega)$ of $u$ is
\begin{align} \label{eq:Su}
S_u(\omega) &= \sum^{\nb}_{i=1} f_u(\omega_i)\delta(\omega-\omega_i) + \sum^{\nb}_{i=1} f_u(-\omega_i)\delta(\omega+\omega_i) \,,
\end{align}
where $f_u(\omega_i)$, $f_u(-\omega_i) > 0$. 
Also by \cite[Prop~1.6.2]{SASTRY89}, 
$S_{\xih}(\omega)= T(-j\omega)S_u(\omega)T\T(j\omega)$, where $T(s)$ is defined in \eqref{eq:xihreachable}. 
Recall that $R[\xih] = \frac{1}{2\pi} \int^{+\infty}_{-\infty} S_{\xih}(\omega) d\omega$ \cite[p.~256]{IOANNOU12}. Combined with \eqref{eq:Su}, we obtain 
\begin{align*}
R[\xih] = &\frac{1}{2\pi}\sum^{\nb}_{i=1}f_u(\omega_i)T(-j\omega_i) T\T(j\omega_i) \\
          &+ \frac{1}{2\pi}\sum^{\nb}_{i=1}f_u(-\omega_i)T(j\omega_i) T\T(-j\omega_i)\,.
\end{align*}
Suppose $R[\xih]$ is not positive definite. Then there exists $v \in \RR^{2\nb}$, $v \neq 0$, such that
\begin{align}
\nonumber
v\T R[\xih]v 
&= v\T \left(\frac{1}{2\pi}\sum^{\nb}_{i=1}f_u(\omega_i)T(-j\omega_i) T\T(j\omega_i) 
       \right. \\
\label{eq:Rxihcond}
&+ \left. \frac{1}{2\pi}\sum^{\nb}_{i=1}f_u(-\omega_i)T(j\omega_i) T\T(-j\omega_i) 
\right) v = 0 \,.
\end{align}
Since $f_u(\omega_i)$, $f_u(-\omega_i)>0$, and each term in the summation is nonnegative (because each term is a complex vector times its conjugate transpose), then \eqref{eq:Rxihcond} implies for $i = 1 \ldots \nb$,
$v\T T(-j\omega_i)T\T (j\omega_i) v = 0$ and 
$v\T T(j\omega_i)T\T(-j\omega_i )v   = 0$.  
This is equivalent to
$v\T T(-j\omega_i) = 0$ and $v\T T(j\omega_i)  = 0$, for $i = 1 \ldots \nb $. In more compact form, we have $v\T \cT = 0$, where 
\begin{align}
\cT &= 
\begin{bmatrix}
 T(s_1) & \ldots & T(s_{\nb}) & T(-s_1) & \ldots & T(-s_{\nb})
\end{bmatrix} \,,
\end{align}
and $s_i = j\omega_i \in \CC$. 
Define $\ol{s} := (s_1,\ldots,s_{\nb})$ and $-\ol{s} := (-s_1,\ldots,-s_{\nb})$. Using \eqref{eq:laplaceexpansion}, we have $\cT = L D$, where 
\begin{align*}
D &= \diag \left( \frac{1}{\al(s_1)f(s_1)}, \ldots,  
                  \frac{1}{\al(s_{\nb})f(s_{\nb})}, \right. \\
  &\qquad \qquad  \left. \frac{1}{\al(-s_{1}) f(-s_{1})}, \ldots,
                  \frac{1}{\al(-s_{\nb})f(-s_{\nb})} \right) 
\end{align*}
\begin{align*}
L &= 
\begin{bmatrix} 
L_1(\ol{s}) & L_1(-\ol{s}) \\
L_2(\ol{s}) & L_2(-\ol{s}) 
\end{bmatrix} \\
L_1(\ol{s}) &=  
\begin{bmatrix}
1   \cdot \beta(s_{1}) & \ldots & 1 \cdot  \beta(s_{\nb}) \\
s_1 \cdot \beta(s_{1}) &  \ldots & (s_{\nb}) \cdot \beta(s_{\nb}) \\
\vdots & & \vdots \\
(s_{1})^{\nb-1} \cdot \beta(s_{1}) & \ldots & 
(s_{\nb})^{\nb-1}\cdot \beta(s_{\nb}) 
\end{bmatrix} \\
L_2(\ol{s}) &= 
\begin{bmatrix}
1 \cdot\alpha(s_1)& \ldots & 1 \cdot \alpha(s_{\nb}) \\
(s_1) \cdot \alpha(s_1) & \ldots & (s_{\nb}) \cdot \alpha(s_{\nb}) \\
\vdots & & \vdots \\
(s_1)^{\nb-1} \cdot \alpha(s_1) & \ldots 
 & (s_{\nb})^{\nb - 1} \cdot \alpha(s_{\nb}) 
\end{bmatrix} \,.
\end{align*}
The polynomials $\al(s)$ and $f(s)$ are Hurwitz, so matrix $D$ is diagonal with nonzero diagonal elements. Hence, $\rnk(\cT) = \rnk(L)$.

Now suppose there exists a vector 
$v = \begin{bmatrix} v_1\T & v_2\T \end{bmatrix}\T \in \RR^{2\nb}$, 
$v \neq 0$ and $v_1$, $v_2 \in \RR^{\nb}$ such that 
$v\T L = \begin{bmatrix} v_1\T & v_2\T \end{bmatrix} L = 0$. 
Equivalently, $L$ is not full rank if there exist polynomials $v_1(s)$, $v_2(s) \neq 0$ with $\deg(v_1(s)), \deg(v_2(s)) \leq \nb - 1$ whose coefficients correspond to the elements of $v_1$, $ v_2$, such that
\begin{align*}
     v_1(s_i) \beta(s_i) + v_2(s_i)\al(s_i) &= 0\,, \\
     v_1(-s_i) \beta(-s_i) + v_2(-s_i)\al(-s_i) &= 0\,, ~~ i = 1 \ldots \nb \,.
\end{align*}
This implies
\begin{align}
\nonumber
v_1(s)\beta(s)+ v_2(s) \al(s) 
&= v_3(s) \prod_{i=1}^{\nb} (s^2 + \omega^2_i) \\
\label{eq:Elliot80}
&= v_3(s) w(s) \,,
\end{align}
for some polynomial $v_3(s)$ of degree at most $n - \nb - 1$. Suppose by way of contradiction, $v_3(s) \equiv 0$. Then we may rewrite \eqref{eq:Elliot80} as $v_1(s)\beta(s) = - v_2(s)\al(s)$.
Since $\al(s)$ and $\beta(s)$ are coprime, $v_1(s)$ must contain the zeros of $\al(s)$. But $\deg(v_1(s)) \leq \nb - 1 < \deg(\al(s))$, a contradiction. Thus, $v_3(s)$ is not identically zero.

The next steps utilizing Sylvester matrices are inspired by arguments in \cite[Lemma~2]{ELLIOTT80}. We rewrite \eqref{eq:Elliot80} as
\begin{align}\label{eq:Elliot80_1}
    \begin{bmatrix}
        v_1\T & v_2\T & -v_3\T
    \end{bmatrix}\begin{bmatrix}
        1 \cdot \beta(s)\\
        s \cdot \beta(s)\\
        \vdots\\
        s^{\nb-1} \cdot \beta(s)\\
        1 \cdot \al(s)\\
        s \cdot \al(s)\\
        \vdots\\
        s^{\nb -1} \cdot \al(s)\\
        1 \cdot w(s)\\
        s \cdot w(s)\\
        \vdots\\
        s^{n - \nb - 1} \cdot w(s)
    \end{bmatrix}&=0\,,
\end{align}
where $v_3$ is a vector in $\RR^{n - \nb }$ containing the coefficients forming $v_3(s)$. There are two cases: (i) $\nb < n < 2\nb$; (ii) $2 \nb \le n$.
We consider only case (i), since the other case is proved analogously. Recalling that $v_1$, $v_2 \in \RR^{\nb}$, and $v_3\in \RR^{n - \nb}$, we can expand \eqref{eq:Elliot80_1} as 
\begin{align}\label{eq:LindCond}
\begin{bmatrix}
v_1\T & v_2\T & -v_3\T
\end{bmatrix} M S(s) &= 0 \,,
\end{align}
where $S(s) := \col( 1, s, \ldots, s^{n + \nb - 1} )$ and 
\begin{align*}
M &= \col( M_1, M_2, M_3 ) \in \RR^{(\nb + n) \times (\nb + n)} \\
M_{1} &=
\left[
\begin{array}{llllllll}
\beta_0 & \ldots & \beta_{n-1} & 0 & \ldots & & 0 \\
\vdots & & \ddots  & & \ddots & & \\
0 & \ldots & \ldots & \beta_0  & \ldots & \beta_{n-1} & 0  \\
\end{array} 
\right] \\
M_{2} &=
\left[
\begin{array}{llllllll}
\al_0 & \ldots & \al_{n-1} & 1 & \ldots & & 0 \\
\vdots & & \ddots  & & \ddots & & \\
0 & \ldots &  \ldots & \al_0 & \ldots & \al_{n-1} & 1 \\
\end{array} 
\right] \\
M_{3} &=
\left[
\begin{array}{llllllll}
w_0 & \ldots & \ldots & \ldots & 1 & 0 & \ldots & 0 \\
\vdots & \ddots & & & \ddots & & \\
0 & \ldots & w_0 & \ldots & \ldots &  & w_{2\nb-1}& 1 \\
\end{array} 
\right]\,.
\end{align*}
The elements of $S(s)$ form a basis for the ring of polynomials of maximum degree $\nb + n - 1$. Thus, \eqref{eq:LindCond} holds iff 
\begin{align*}
\begin{bmatrix}
v_1\T & v_2\T & -v_3\T
\end{bmatrix} M &= 0\,.
\end{align*}
Now we argue that matrix $M$ is full rank for almost all choices of frequencies $\omega_1 \ldots \omega_{\nb}$. Notice that $M$ is a partial \textit{generalized Sylvester resultant matrix} for polynomials $\al(s)$, $\beta(s)$, and $w(s)$. 
By assumption, the polynomials $\al(s)$ and $\beta(s)$ are coprime. Furthermore, $\al(s)$ is Hurwitz and $w(s)$ has purely complex roots, which makes them coprime. Next, the set $\{ s \in \CC ~|~ \beta(s) = 0 \}$ has measure zero, while the roots of polynomial $w(s)$ correspond to $\nb$ distinct, randomly chosen frequencies $\omega_i$ as well as their negative counterpart. Thus, for almost every choice of frequencies $\omega_i$, $i = 1\ldots \nb$, we have that $\beta(j\omega_i) \neq 0$, which implies $w(s)$ and $\beta(s)$ are coprime. We conclude $\alpha(s)$, $\beta(s)$, and $w(s)$ have no common roots, so by Barnett's Theorem \cite[p.~39]{BARNETT83}, $M$ is nonsingular. This means there does not exist a non-zero $v \in \RR^{2\nb}$ such that $v\T\cT = 0$, implying $R[\xih]$ is positive definite, and $\xih$ is PE.

\subsection{Proof of Theorem~\ref{thm:params}}

First we study stability of the subsystem \eqref{eq:thtdot} and \eqref{eq:xitdot}. Consider the Lyapunov function 
\begin{align}
V &= \frac{1}{2} \tht\T \tht + \gamma\0 \xit\T P \xit \,,
\end{align}
where $\gamma\0 > 0$ is to be determined, and $P$ satisfies the Lyapunov equation $\Fbz\T P + P \Fbz = -I$. Computing the Lie derivative,
\begin{align*}
\dot{V} 
& \leq  
-\frac{\gamma}{2} \| \tht\T \xih \|^2 
-\biggl( \gamma\0 - \frac{\gamma}{2} \| \theta \|^2 \biggr) \| \xit \|^2 \,,
\end{align*}
where we have applied Young's inequality and Cauchy-Schwarz.
Choosing $\gamma\0 > \gamma \| \theta\|^2$, we obtain $ \dot{V} \leq 0$.
By \cite[Theorem~A.5]{KRSTIC95}, the equilibrium $(\tht,\xit) = (0,0)$ is globally uniformly stable (GUS).

Second, we show $\xis(t)$ is PE. Using $\nb = n$, \eqref{eq:FG} and \eqref{eq:nominalparams}, the model \eqref{eq:overp2} becomes
\begin{align}
\label{eq:overp5}
\dot{\xi} 
&= \begin{bmatrix}
    A & G \beta\T \\
    0 & F
\end{bmatrix} \xi + \begin{bmatrix} 0 \\ G \end{bmatrix} u 
= \Fb(\theta) \xi + \Gb u \,.
\end{align}
where $A = A\0 + B\0 (\alpha)\T$.
Since both $A$ and $F$ are Hurwitz, then $\Fb(\theta)$ is Hurwitz. By Theorem~\ref{thm:order}, the controllable subspace of \eqref{eq:overp5} satisfies $\dim(\cC) = 2 n$, so the pair $(\Fb(\theta),\Gb)$ is controllable. Since $u(t)$ is piecewise continuous, bounded, and SR of order $2n$, we can apply \cite[Theorem~6.3]{NARENDRA89} to conclude that $\xis(t)$ is PE. 

Next consider the nominal error model
\begin{align}
\label{eq:nominal}
\dot{\tht} &= - \gamma \left( \xis(t) \xis\T(t) \right) \tht \,.
\end{align}
Since $\xis(t)$ is stationary and PE, by \cite[Theorem~2.5.1]{SASTRY89}, the equilibrium $\tht = 0$ is globally exponentially stable (GES). Noting the Jacobian of \eqref{eq:nominal} is bounded (since \eqref{eq:nominal} is globally Lipschitz), by a global version of \cite[Theorem~4.14]{KHALIL02}, there exists $V_0:  [0,\infty) \times \RR^{2n}\rightarrow\RR$ and 
$\al_i > 0$ such that
\begin{subequations}
\begin{align}
\al_1\|\tht\|^2 \leq V_0(t,\tht) 
&\leq \al_2 \|\tht\|^2\\
\frac{\p V_0}{\p t} + \frac{\p V_0}{\p \tht} (- \gamma  \xis(t) \xis\T(t) \tht) 
&\leq - \al_3 \|\tht\|^2 \\
\left\|\frac{\partial V_0}{\p \tht} \right\| 
&\leq \al_4 \|\tht\| \,.
\end{align}    
\end{subequations}

Fix $r_0 > 0$ and let $(\tht(t_0), \Sigt(t_0), \xit(t_0)) \in \cB(r_0)$. Consider the Lyapunov function
\begin{align*}
V_1(t,\tht,\xit) &:= V_0(t,\tht) + \gamma_1 \xit\T P \xit \\
V_2(\Sigt) &= \| \Sigt \|^2_F \\
V(\cdot) &= V_1(t,\tht,\xit) + V_2(\Sigt) \,,
\end{align*} 
where $\gamma_1 > 0$ is to be determined, and $\| \cdot \|_F$ denotes the Frobenius norm. By the GUS argument for $(\tht,\xit) = (0,0)$, there exists $r(r_0)$ such that 
\begin{align*}
\| \tht(t) \| &\le r(r_0) \,, \quad \| \xit(t) \| \le r(r_0) \,,
\quad t \ge t_0 \ge 0 \,.
\end{align*}
Also, define $\cs := \| \xis \|_{\cL_\infty}$. Using properties of the Frobenius norm (triangle inequality, sub-multiplicative), we have
\begin{align*}
\| g(t,\xit) \| &\le \| g(t,\xit) \|_F
\leq \| \xit \xis\T(t) \|_F + \| \xis(t) \xit\T \|_F + \| \xit \xit\T \|_F \\
&\leq 2 \cs \|\xit \|+ \|\xit\|^2 
 \leq ( 2 \cs + r(r_0) ) \| \xit \| =: c_g(r_0) \| \xit \| \,.
\end{align*} 
Now we compute the Lie derivative of $V_1(\cdot)$: 
\begin{align*}
\dot{V}_1(\cdot) 
&\le - \al_3 \| \tht \|^2 - \gamma_1 \| \xit \|^2 
 + \alpha_4 \gamma  \| g(t,\xit) \| \| \tht \|^2 \\
&+ \alpha_4 \gamma \| \theta \| \| \tht \|  \| \xit \|^2
 + \alpha_4 \gamma \| \theta \| \cs \| \tht \| \| \xit \| \\
&\le - \al_3 \| \tht \|^2 - \gamma_1 \| \xit \|^2 \\
&+ \alpha_4 \gamma 
\bigl( c_g(r_0) r(r_0) + r(r_0) \| \theta \| + \cs \| \theta \| \bigr)
\| \tht \| \| \xit \| \\
&=: - \al_3 \| \tht \|^2 - \gamma_1 \| \xit \|^2 
+ c_1(r_0) \| \tht \| \| \xit \| \,.
\end{align*}
Apply the Peter Paul inequality:
$\| \tht \| \| \xit \| \leq \frac{\|\tht\|^2}{2 \nu} +\frac{\nu \| \xit \|^2}{2}$, with $\nu > 0$. Then 
\begin{align*}
\dot{V}_1(\cdot) &\le 
- \left( \al_3 -\frac{c_1(r_0)}{2\nu} \right) \|\tht\|^2 
- \left( \gamma_1 - \frac{c_1(r_0) \nu}{2} \right) \| \xit \|^2 \,.
\end{align*}
Choose $\nu = \frac{c_1(r_0)}{\alpha_3} > 0$. Then we get
\begin{align*}
\dot{V}_1(\cdot) &\le 
- \frac{\al_3}{2} \|\tht\|^2 
- \left( \gamma_1 - \frac{c_1^2(r_0)}{2 \alpha_3} \right) \| \xit \|^2 \,.
\end{align*}
Next compute the Lie derivative of $V_2(\Sigt)$:
\begin{align*}
\dot{V}_2(\Sigt)
&= - 2 \veps  \| \Sigt \|^2_F 
+ 2 \veps \tr( \Sigt\T g(t,\xit) ) \\
&\leq - 2 \veps  \| \Sigt \|^2_F 
+ 2 \veps  \| \Sigt \|_F \| g(t,\xit) \|_F \\
&\leq - \veps  \| \Sigt \|^2_F 
+ \veps \| g(t,\xit) \|_F^2 \\
&\leq - \veps  \| \Sigt \|^2_F + \veps c_g^2(r_0) \| \xit \|^2
\end{align*}
where we used Cauchy-Schwarz inequality and Young's inequality.
Returning to $V(\cdot)$, define the constants
\begin{align*}
\gamma_1(r_0) 
&:= 1 + \frac{c_1^2(r_0)}{2 \alpha_3} + \veps c_g^2(r_0) \\
\lambda 
&:= \min \{ \frac{\alpha_3}{2} , 1, \veps \} \\
\beta_1(r_0) 
&:= \min \{ \alpha_1, 1, \gamma_1(r_0) \lambda_{\min}(P) \} \\
\beta_2(r_0) 
&:= \max \{ \alpha_2, 1, \gamma_1(r_0) \lambda_{\max}(P) \} \,.
\end{align*}
Then we verify
\begin{align*}
\beta_1(r_0) \| ( \tht, \vect(\Sigt),  \xit ) \|^2 
&\le V(\cdot) \le \beta_2(r_0) \| ( \tht, \vect(\Sigt),  \xit ) \|^2 
\end{align*}
where $\vect(\Sigt)$ is the vectorization of $\Sigt$.
Also, we compute
\begin{align*}
    \dot{V} & \leq - \frac{\al_3}{2} \|\tht\|^2 
- \left( \gamma_1(r_0) - \frac{c_1^2(r_0)}{2 \alpha_3} - \veps c_g^2(r_0) \right) \| \xit \|^2  \\
            &- \veps  \| \Sigt \|^2_F \\
            & \leq - \frac{\al_3}{2} \|\tht\|^2  -  \| \xit \|^2 - \veps  \| \Sigt \|^2_F\\
            & \leq - \lambda \| ( \tht, \vect(\Sigt),  \xit ) \|^2\\
            & \leq - \frac{\lambda}{\beta_2 (r_0)} V \,.
\end{align*}
An application of the Comparison Lemma \cite[Lemma~3.4]{KHALIL02} completes the proof.

\section{Proof of the Laplace Expansion in Theorem~\ref{thm:order}} 
\label{sec:Laplace}

We apply the Laplace expansion along the first row of \eqref{eq:block3} to compute its determinant
\begin{align}
\label{eq:expansion}
\det(M_{22}) 
&= \sum^{\bar{n}}_{j=1} A_{1j}C_{1j}
 = \sum^{\bar{n}}_{j=1} (-1)^{(j+1)}\bar{\beta}_{(j-1)} \det (m_{1j})
\end{align}
where $A_{1j} = \bar{\beta}_{(j-1)}$ is the $j^{\text{th}}$ entry of the first row of $M_{22}$ and $C_{1j} = (-1)^{(j+1)}\det(m_{1j})$ is the corresponding cofactor.
The submatrix $m_{11}$ has a lower diagonal structure of the form
\begin{equation}
m_{11} =  
 \begin{bmatrix}
 1 & 0 & 0 & 0 & \ldots & 0\\
        -\lambda & 1 & 0 & 0 & \ldots & 0\\
        0 &  -\lambda & 1 & 0  & \ldots & 0\\
        \vdots & &\ddots & \ddots \\
        0 & &&-\lambda & 1 & 0\\
        0 & &&0&-\lambda & 1 
 \end{bmatrix} \,,
\end{equation}
then $\det(m_{11}) = 1$. Notice that the submatrices $m_{1j} \in \RR^{(\nb - 1) \times (\nb - 1)}$, with $j\geq 2$, have a special block diagonal structure of the form

\begin{align*}
m_{1j} & =   
       \begin{bmatrix}
            B_{1j} & 0 \\
        0 & B_{2j}
       \end{bmatrix} \\
   &= 
\left[
\begin{array}{ccccc|ccccccc}
           -\lambda & 1 & 0 & 0 &&&&&\ldots & 0\\
        0 & -\lambda & 1 &  0 &&&&& \ldots & 0\\
        \vdots & &\ddots & \ddots &&&\\
        0 & &  & -\lambda & 1 &  0 &&& \ldots & 0\\
        0 & & & & -\lambda & 0 && &\ldots & 0 \\
        \hline
        0 & &&&0 & 1 & 0 & &\ldots & 0 \\
        0 & &  & &0&-\lambda & 1 &  0 & \ldots & 0\\
        \vdots && & & & &\ddots & \ddots \\
               0 & &  &&&&0&-\lambda & 1&0\\
        0 & &  & &&&&0&-\lambda & 1\\
\end{array}
\right]\,,    
\end{align*}
where $B_{1j} \in \RR^{(j-1) \times (j-1)}$ and $B_{2j} \in \RR^{(\nb- j) \times (\nb - j)}$ are upper and lower triangular matrices, respectively. Hence,
\begin{align}\label{eq:minor_f}
   \det(m_{1j}) = \det(B_{1j}) \det(B_{2j}) = (-1)^{j-1} \lambda^{j-1} \times 1 \,.
\end{align}
Substituting \eqref{eq:minor_f} into \eqref{eq:expansion} yields
\begin{align*}
 \det(M_{22})  = \sum^{\bar{n}}_{j=1}\bar{\beta}_{(j-1)}\lambda^{j-1} =  \beta_e(\lambda)\,.
\end{align*}

\section{Conclusion} 

The paper has taken a step toward explaining in control theoretic terms how the brain may be able to learn low-order LTI models. Our framework utilizes stable linear filters, a gradient adaptation law, and an SVD computation, all biologically plausible computations in the brain. The model order identification algorithm is a proxy for {\em episodic learning} over successive tasks, as well as {\em Bernstein's principle}, that learning progresses from simple low order behaviors to higher order behaviors \cite{NEWELL01}. The variable $\veps$ plays an important role in setting the timescale over which learning takes place, acknowledging that the brain utilizes a vast separation of timescales \cite{KIEBEL08}.
The framework focuses on learning low order SISO models of LTI systems, a capability most closely associated with the oldest motor systems phylogenetically. A significant extension of the framework is to MIMO nonlinear systems so that we can study how the brain learns nonlinear models for regulation of arm movement, locomotion, and manipulation of tools. 

\bibliographystyle{plain}
\bibliography{SIDBib}

\end{document}